\documentclass[review,number]{elsarticle}

\usepackage{amsmath,amssymb}
\usepackage{commath} 
\usepackage{amsxtra} 
\usepackage{booktabs}
\usepackage{subfig}

\usepackage{colonequals}

\usepackage{array}

\usepackage[noprefix]{nomencl}
\makenomenclature
\renewcommand{\nomname}{NOTATION}
\makeatletter
\renewenvironment{thenomenclature}{%
  \@ifundefined{chapter}%
  {
    \section*{\nomname}
    \if@intoc\addcontentsline{toc}{section}{\nomname}\fi%
  }%
  {
    \chapter*{\nomname}
    \if@intoc\addcontentsline{toc}{chapter}{\nomname}\fi%
  }%

  \nompreamble
  \list{}{%
    \labelwidth\nomlabelwidth
    \leftmargin\labelwidth
    \advance\leftmargin\labelsep
    \itemsep\nomitemsep
    }}
\makeatother

\usepackage{algorithm}
\usepackage{algpseudocode}
\floatname{algorithm}{Procedure}

\newcommand{\subeq}[2]{\eqref{#1}$_{#2}$}

\newcommand{\argmin}{\operatornamewithlimits{argmin}}
\newcommand{\sbt}{\text{subject to }}

\newcommand{\half}{\frac{1}{2}}
\newcommand{\thalf}{\tfrac{1}{2}}

\newcommand{\NDOF}{N_{\text{DOF}}}
\newcommand{\Nsig}{N_\sigma}
\newcommand{\Nzeta}{N_\zeta}
\renewcommand{\NG}{N_\text{G}}
\newcommand{\Ny}{N_\text{y}}

\newcommand{\psic}{\psi^\text{c}}
\newcommand{\phic}{\phi^\text{c}}
\newcommand{\Pic}{\Pi^\text{c}}

\newcommand{\un}{u^n}
\newcommand{\unp}{u^{n+1}}
\newcommand{\vn}{v^n}\newcommand{\vnp}{v^{n+1}}
\newcommand{\sigman}{\sigma^n}
\newcommand{\sigmanp}{\sigma^{n+1}}
\newcommand{\zetan}{\zeta^n}
\newcommand{\zetanp}{\zeta^{n+1}}
\newcommand{\pnp}{p^{n+1}}

\newcommand{\lambdabar}{\bar{\lambda}}

\newcommand{\lambdabarnp}{\lambdabar^{n+1}}

\newcommand{\Dt}{\Delta t}

\newcommand{\bsigma}{b_\sigma}
\newcommand{\bzeta}{b_\zeta}
\newcommand{\muprsc}{\mu^\text{prsc}}
\newcommand{\Bprsc}{B^\text{prsc}}

\newcommand{\sigmam}{\sigma_m}
\newcommand{\zetam}{\zeta_m}
\newcommand{\lambdam}{\lambda_m}
\newcommand{\bsigmam}{b_{\sigma m}}
\newcommand{\bzetam}{b_{\zeta m}}

\newcommand{\sumM}{\sum^{\NG}_{m=1}}
\newcommand{\forM}{m=1,\ldots,\NG}

\newcommand{\calL}{\mathcal{L}}

\newcommand{\lamtilm}{\tilde{\lambda}_m}
\newcommand{\sigtildem}{\tilde{\sigma}_m}
\newcommand{\zetatildem}{\tilde{\zeta}_m}

\newcommand{\sigzetavec}{\del{\begin{smallmatrix} \sigma \\ \zeta \end{smallmatrix}}}
\newcommand{\mulamvec}{\del{\begin{smallmatrix} \mu \\ \lambda \end{smallmatrix}}}
\newcommand{\sigzetavecm}{\del{\begin{smallmatrix} \sigma_m \\ \zeta_m \end{smallmatrix}}}
\newcommand{\mulamvecm}{\del{\begin{smallmatrix} \mu \\ \lambda_m \end{smallmatrix}}}

\newcommand{\LTosigmam}{\pd{\calL}{\sigmam}}
\newcommand{\LTozetam}{\pd{\calL}{\zetam}}

\newcommand{\sigmstar}{\sigma^*_m}
\newcommand{\zetamstar}{\zeta^*_m}

\newcommand{\dmu}{\delta\mu}

\newcommand{\Kss}{K_{\sigma\sigma}}
\newcommand{\Ksx}{K_{\sigma\zeta}}

\newcommand{\Kxx}{K_{\zeta\zeta}}
\newcommand{\Phis}{\Phi_\sigma}
\newcommand{\Phix}{\Phi_\zeta}
\newcommand{\Kbarxx}{\bar{K}_{\zeta\zeta}}
\newcommand{\Phibars}{\bar{\Phi}_\sigma}
\newcommand{\Phibarx}{\bar{\Phi}_\zeta}

\newcommand{\Kbar}{\bar{K}}

\newcommand{\calI}{\mathcal{I}}

\newcommand{\iter}{\mathtt{iter}}
\newcommand{\maxiter}{\mathtt{MAXITER}}

\DeclareMathAlphabet{\mathpzc}{OT1}{pzc}{m}{it}

\newcommand{\nback}{\mathtt{nback}}
\newcommand{\maxbacktrack}{\mathtt{MAXBACKTRACK}}

\newcommand{\sigmay}{\sigma_\text{y}}

\newcommand{\sqrttt}{\sqrt{\frac{2}{3}}}

\newcommand{\zetakh}{\zeta_\text{kh}}
\newcommand{\zetaih}{\zeta_\text{ih}}
\newcommand{\psicih}{\psic_\text{ih}}
\newcommand{\bzetakh}{b_{\zeta_\text{kh}}}
\newcommand{\bzetaih}{b_{\zeta_\text{ih}}}
\newcommand{\sigmatr}{\sigma^\text{trial}}
\newcommand{\zetatr}{\zeta^\text{trial}}
\newcommand{\zetakhtr}{\zeta^\text{trial}_\text{kh}}
\newcommand{\zetaihtr}{\zeta^\text{trial}_\text{ih}}
\newcommand{\ys}{y_\sigma}
\newcommand{\ykh}{y_{\zetakh}}
\newcommand{\LambdaCHP}{\Lambda_{(C+H)P}}
\newcommand{\psiih}{\psi_\text{ih}}
\newcommand{\alphih}{\alpha_\text{ih}}

\biboptions{sort&compress}

\begin{document}

\begin{frontmatter}

    \title{An optimization-based reformulation of the classical displacement approach for state update of non-linear material models}

    \author[UB]{Z.~Lotfian}
    \ead{zahrasad@buffalo.edu}

    \author[UB]{M.~V.~Sivaselvan\corref{cor1}}
    \ead{mvs@buffalo.edu}

    \cortext[cor1]{Corresponding author}
    \address[UB]{Department of Civil, Structural and Environmental Engineering, University at Buffalo, 212 Ketter Hall, Buffalo, NY 14260, U.S.A.}

    \begin{abstract}
        In this paper, we build on recent work using a mathematical programming approach for incremental state update in analysis of non-linear mechanics models.
        In particular, we consider quasi-static analysis of continuum problems in the linearized kinematics regime,
        with non-linear material models described using \emph{convex} energy functions.
        We find in this case that the classical displacement-based nested approach for incremental state update
        can be reformulated as solving a reduced \emph{dual} optimization problem.
        This reformulation provides insights into the working of the algorithm, and eliminates the need for some heuristics.
        An important purpose of this paper is to further illustrate the unifying nature of the mathematical programming approach.
        We therefore present relationships with several of these types of algorithms recently presented in the literature for incremental state update.
    \end{abstract}


\end{frontmatter}

\section{Introduction}\label{sec:Intro}
Analysis of mechanics models with complex non-linear material behavior arises in different applications such as seismic response simulation of structures.
In a recent paper \cite{SivaselvanNLM2013}, we discussed the manner in which when
such material behavior is described using an energy approach \cite{Halphen1975,HoulsbyPurzin,ErlicherEndochronic2006},
\emph{state update} of the model in each increment of a non-linear analysis can be cast into various \emph{mathematical programs}.
Depending on the type of material model, the corresponding mathematical program could assume different forms
such as convex optimization, complementarity problem etc.
This mathematical programming approach represents an alternative to the classical displacement-based approach.
In the latter, state update is carried out in a nested fashion with displacements computed at a global level,
and stresses and other internal variables computed at a material-point level.
It was also seen in \cite{SivaselvanNLM2013} that models of non-linear behavior that arose historically from disparate contexts,
for example the Preisach model \cite{Preisach1935,SpanosSteadyStatePreisach} and
the Bouc-Wen model \cite{Bouc1971,Wen1976,BursiBoucWen2008,BoucWenSurvey},
could be interpreted within the mathematical programming framework.

In the present paper, we further explore the unifying nature of the mathematical programming approach.
We consider quasi-static continuum problems in the linearized kinematics regime, with non-linear material behavior described by convex energy functions.
In this case, we find that the classical displacement-based nested approach can be reformulated as solving a reduced \emph{dual} optimization problem.
Some differentiability results related to convex optimization problems play an important role in developing this reformulation.
This reformulation in turn provides insights into the working of the algorithm by means of some geometric constructions such as Figure \ref{fig:linesearch},
and obviates some heuristics that are otherwise used.
We envision that such insights will guide development of algorithms for more complex models,
for example via successive convex programs \cite{AnitescuSeqConvex,Garcea2011,MilaniMasonry2011,KrabbenhoftNonassociated}.

A particular case of non-linear material models with convex energy functions is rate-independent plasticity models.
For this case, incremental state update at the material point level is known as return mapping, and has a well-known optimization format \cite{SimoHughesCompInelasticity}.
However, the optimization structure at the global level is not commonly recognized or utilized.
The purpose of this paper is two-fold --- (1) to reformulate the classical displacement-based nested approach for incremental state update
as an optimization problem (for quasi-static problems with convex energy functions), and (2) identify relationships with several other algorithms recently described in the literature.

This paper is organized as follows.
First in section \ref{sec:GoverningEq}, governing equations presented in \cite{SivaselvanNLM2013} for single-degree-of-freedom systems with energy-based material models
are extended to general spatially discretized systems. With some restrictions, when discretized in time, these equations lead to an optimization problem in each time increments.
In section \ref{sec:Opt}, three forms of this optimization problem are discussed --- a primal problem, a dual problem, and a reduced dual problem.
The reduced dual problem is a reformulation of the classical displacement-based nested approach,
and consists of solving an optimization sub-problem at the integration point level.
Solution of the reduced dual problem is the subject of section \ref{sec:reformulate}.
Differentiability of the objective function is considered and derivatives are obtained,
so that Newton's method for unconstrained minimization can be applied as described in section \ref{sec:reformulate:algorithm}.
In section \ref{sec:return-map}, some special classes of material models are discussed, where the integration point level optimization problem admits a simpler solution process.
A numerical example is presented in section \ref{sec:NumExample} to illustrate the working of the algorithm of section \ref{sec:reformulate:algorithm}.
Lastly in section \ref{sec:relationships}, relationships between various solution strategies presented in recent years based on the mathematical programming approach are discussed.

\section{Governing equations and time discretization}\label{sec:GoverningEq}
\begin{figure}
    \centering
    \subfloat[Indicator function of a convex subset $\mathcal{C}$ of $\mathbb{R}^2$]
    {
        \resizebox{0.4\textwidth}{!}{\includegraphics{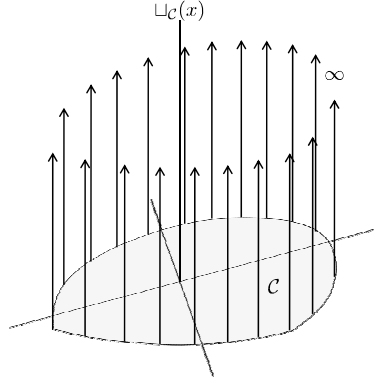}}
        \label{fig:IndFunc:func}
    } \\
    \subfloat[Subdifferential of the indicator function of $\mathcal{C}$;
              The subdifferential is shown at three points --- $x_1$ in the interior if $\mathcal{C}$, $x_2$ where the boundary is smooth, and $x_3$ where the boundary is not smooth]
    {
        \resizebox{0.8\textwidth}{!}{\includegraphics{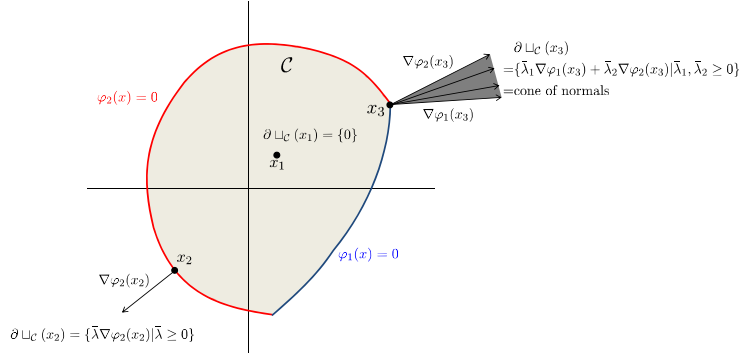}}
        \label{fig:IndFunc:subdifferential}
    }
    \caption{Indicator function and its subdifferential for a convex set $\mathcal{C} = \{x\in\mathbb{R}^n\,|\,\varphi_i(x) \le 0, i=1,2\}$}
    \label{fig:IndFunc}
\end{figure}

We begin by formulating the governing equations for a continuum model that has been spatially discretized, for example by the finite element element method.
We employ a class of non-linear material models described completely by certain energy functions.
As discussed in \cite{SivaselvanNLM2013}, such a representation of material models is based on
the \emph{generalized standard material} framework \cite{Halphen1975,MielkeGSM2006},
and is also closely related to the hyperplasticity framework \cite{CollinsHoulsby,HoulsbyPurzin}.
In this energy framework, a material model is characterized by two convex (possibly nonsmooth) functions
--- a stored energy function $\psi(\epsilon,\alpha)$, and a dissipation function $\phi(\dot{\epsilon},\dot{\alpha})$,
where $\epsilon$ is the strain, and $\alpha$ is a vector of kinematic internal state variables (such as equivalent plastic strain, damage etc.).
We refer to the Fenchel-Legendre transforms of these functions as the complementary stored energy function $\psic(\sigma,\zeta)$,
and the complementary dissipation function $\phic(\sigma,\zeta)$, where $\sigma$ is the stress and the $\zeta$ is the generalized stress conjugate to $\alpha$.
The governing equations discussed in \cite{SivaselvanNLM2013} for simple single-degree-of-freedom dynamic models
can be generalized to spatially discretized continuum models as
\begin{equation}
    \begin{alignedat}{5}
        \text{momentum conservation}&:\quad&
        &\underbrace{\od{}{t}\nabla\del{\thalf \dot{v}^\top M \dot{v}}}_{\text{\tiny{inertia force}}}&
        &{}+{}\underbrace{\nabla\del{\thalf \dot{v}^\top C_d \dot{v}}}_{\text{\tiny{damping force}}}&
        &{}+{}\underbrace{B^\top\sigma}_{\text{\tiny{element force}}}&
        &=p \\
        \text{deformation compatibility}&:\quad&
        &\underbrace{\od{}{t}\partial_1\psic(\sigma,\zeta)}_{\text{\tiny{elastic deformation rate}}}&
        &{}+{}\underbrace{\partial_1\phic(\sigma,\zeta)}_{\text{\tiny{\parbox{0.7in}{plastic deformation rate}}}}&
        &{}-{}\underbrace{B\,v}_{\text{\tiny{\parbox{0.65in}{total deformation rate}}}}&
        &\ni 0 \\
        \text{internal variable evolution}&:\quad&
        &\underbrace{\od{}{t}\partial_2\psic(\sigma,\zeta)}_{\text{\tiny{\parbox{0.7in}{reversible internal variable rate}}}}&
        &{}+{}\underbrace{\partial_2\phic(\sigma,\zeta)}_{\text{\parbox{0.7in}{\tiny{irreversible internal variable rate}}}}&
        &&
        &\ni 0
    \end{alignedat}
    \label{eq:GoverningEqGeneral}
\end{equation}
Here, $v \in \mathbb{R}^{\NDOF}$ is the vector of free velocity components,
$M, C_d \in \mathbb{R}^{\NDOF \times \NDOF}$ are respectively positive semi-definite mass and damping matrices,
and $p \in \mathbb{R}^{\NDOF}$ is the external load vector.
$\sigma \in \mathbb{R}^{\Nsig\NG}$ is the stress vector comprised of $\Nsig$ components at each of the $\NG$ material points.
Similarly, $\zeta \in \mathbb{R}^{\Nzeta\NG}$ is the vector of generalized stress internal variable.
$\psi^c, \phi^c : \mathbb{R}^{\Nsig\NG} \times \mathbb{R}^{\Nzeta\NG} \rightarrow \mathbb{R}$
are the complementary stored energy and complementary dissipation functions.
$B \in \mathbb{R}^{\Nsig\NG\times\NDOF}$ is the linearized deformation-displacement matrix.
$\nabla$ denotes the gradient, and $\partial$ denotes the subdifferential of a nonsmooth convex function \cite{BorweinConvexAnalysis}.
$\nabla_q$ and $\partial_q$ denote gradient and subdifferential with respect to argument $q$ of the function.
Equations \subeq{eq:GoverningEqGeneral}{{2,3}} are inclusions because $\psic$ and $\phic$ could be nonsmooth, and hence have set-valued derivatives.
We make the following remarks on equation \eqref{eq:GoverningEqGeneral}.
\begin{enumerate}
    \item As suggested by their format, equations \eqref{eq:GoverningEqGeneral} can be obtained as
        Euler-Lagrange equations of a generalized Hamilton's principle \cite{SivaselvanNLM2013,Sivaselvan2006}.
    \item The terminology we have used above in reference to spatial discretization follows the common setting,
        where displacements and velocities are associated with nodes, and stresses and other internal variables with element integration points.
        However, the formulation itself is not restricted to this type of spatial discretization.
        For example, in a stress-based finite element, $\sigma$ would denote components of the stress approximation rather than the integration point values,
        or in an assumed-strain finite element, $v$ would include components of the element-level assumed strain coefficients.
    \item For notational convenience, we take that equations \subeq{eq:GoverningEqGeneral}{{2,3}}
        have been multiplied through by the appropriate integration weights and element Jacobian determinants.
        This explains our use of the term ``deformation rate'' instead of ``strain rate'' in equations \subeq{eq:GoverningEqGeneral}{{2,3}},
        and the reason why the element nodal forces can be written simply as $B^\top\sigma$ in equation \subeq{eq:GoverningEqGeneral}{1}.
        Later in section \ref{sec:reformulate:derivatives}, the tangent stiffness matrix is of the form $B^\top \square B$.
    \item In the above, we consider a restricted class of material models for which $\psi$ and $\phi$ are convex.
        One source of non-convexity is so-called \emph{non-associated flow} in plasticity models.
        It can be shown that in some instances of non-associated flow models,
        the complementary dissipation function can be written as $\phic(\sigma,\zeta;\alpha)$
        with explicit dependence on the kinematic internal variables \cite{ErlicherDruckerPrager2005}.
        Equations \eqref{eq:GoverningEqGeneral} then still apply, as discussed for simple single-degree-of-freedom systems in \cite{SivaselvanNLM2013}.
        Explicit dependence on $\alpha$ then suggests solving the incremental state update as a succession of convex subproblems,
        computing $\alpha$ in an outer stage.
    \item When nonlinear kinematics is considered, the term $B$ in equations \subeq{eq:GoverningEqGeneral}{{1,2}}
        is replaced by $\Dif\mathcal{E}(u)$. $\mathcal{E}$ is the nonlinear deformation-displacement map (i.e., the function that maps node displacements to element strains),
        $u$ is the vector of free displacement components, and $\Dif$ is the derivative operator. Thus $B = \Dif\mathcal{E}(0)$.
        In some instances, the kinetic energy in the first term of equation \subeq{eq:GoverningEqGeneral}{1} may also be of the form $\thalf v^\top M(u) v$.
        This again suggests a successive convex programming approach, updating $u$ in an outer stage.
\end{enumerate}
\nomenclature{$v$, $\vn$, $\vnp$}{Velocities at free DOF, and its values at times $n$ and $n+1$ $(\in \mathbb{R}^{\NDOF})$}
\nomenclature{$\sigma$, $\sigman$, $\sigmanp$, $\sigmam$}{Stress at a material point $(\in \mathbb{R}^{\Nsig})$
                                                          or collection over all material points $(\in \mathbb{R}^{\Nsig\NG})$ depending on context;
                                                          superscripts denote time increment index, and subscripts material point index}
\nomenclature{$\epsilon$}{Strain at a material point}
\nomenclature{$\zeta$, $\zetan$, $\zetanp$, $\zetam$}{Generalized stress internal variable at a material point $(\in \mathbb{R}^{\Nzeta})$
                                              or collection over all material points $(\in \mathbb{R}^{\Nzeta\NG})$ depending on context;
                                              superscripts denote time increment index, and subscripts material point index}
\nomenclature{$\alpha$}{Kinematic internal variable at a material point}
\nomenclature{$\psic$}{Complementary stored energy function for a material point, or for entire model, depending on context}
\nomenclature{$\phic$}{Complementary dissipation function for a material point, or for entire model, depending on context}
\nomenclature{$\psi$}{Stored energy function for a material point}
\nomenclature{$\phi$}{Dissipation function for a material point}
\nomenclature{$B$}{Linearized strain-displacement matrix $(\in \mathbb{R}^{\Nsig\NG\times\NDOF})$}
\nomenclature{$M$}{Mass matrix $(\in \mathbb{R}^{\NDOF \times \NDOF})$}
\nomenclature[$Cd$]{$C_d$}{Damping matrix $(\in \mathbb{R}^{\NDOF \times \NDOF})$}
\nomenclature{$p$, $\pnp$}{External load vector, and its values at time $n+1$ $(\in \mathbb{R}^{\NDOF})$}
\nomenclature[$NDOF$]{$\NDOF$}{Number of free degrees of freedom}
\nomenclature[$N\sigma$]{$\Nsig$}{Number of stress components per material point (for example, 3 in 2D problems and 6 in 3D problems)}
\nomenclature[$N\zeta$]{$\Nzeta$}{Number of internal variables per material point (for example, 4 for 2D elastoplasticity with one yield condition and combined isotropic and kinematic hardening)}
\nomenclature[$NG$]{$\NG$}{Number of integration points}
\nomenclature[$E$]{$\mathcal{E}$}{Nonlinear deformation-displacement map used only in the remarks following equation \eqref{eq:GoverningEqGeneral}}
\nomenclature{$u$, $\un$, $\unp$}{Displacement at free DOF, and its values at times $n$ and $n+1$ $(\in \mathbb{R}^{\NDOF})$}
\nomenclature[$D$]{$\Dif$, $\Dif_q$}{Derivative of a differentiable function; subscript $q$ denotes with respect to argument $q$}
\nomenclature[$~$]{$\nabla$, $\nabla_q$}{Gradient of a differentiable function; subscript $q$ denotes with respect to argument $q$}
\nomenclature[$~$]{$\partial$, $\partial_q$}{Subdifferential of a nonsmooth convex function; subscript $q$ denotes with respect to argument $q$}
\nomenclature{$t$}{Time}
For the sake of concreteness, in what follows, we introduce the following additional restrictions.
\begin{enumerate}
    \item The complementary stored energy function $\psic$ is smooth.
        Nonsmooth stored energy functions, such as those resulting from approximating tension- or compression-only behavior
        and some forms of damage, can however be treated in a similar manner to nonsmooth dissipation functions
        as described below \cite{SivaselvanNLM2013}.
    \item The complementary dissipation function is of the form
        \begin{subequations} 
            \label{eq:YieldFunc}
            \begin{equation}
                \phic(\sigma,\zeta) = \sqcup_\mathcal{C}(\sigma,\zeta)
                \label{eq:YieldFunc:indfunc}
            \end{equation}
            where
            $\mathcal{C}$ is a convex set given by
            \begin{equation}
                \mathcal{C} = \{ (\sigma,\zeta) | \varphi(\sigma,\zeta) \le 0\}
                \label{eq:YieldFunc:set}
            \end{equation}
            $\sqcup_\mathcal{C}$ is its indicator function \cite{BorweinConvexAnalysis} (Figure ,
            and $\varphi$ is a smooth vector-valued function, each component of which is convex.
            In plasticity models, $\varphi$ is the yield function, and $\mathcal{C}$ is the elastic region.
        \end{subequations} 
\end{enumerate}
We note that this restricted framework is sufficient to describe nonlinear elastoplasticity.
We assume for simplicity of presentation that there is an identical number $\Ny$ of yield functions at each material point.
Thus $\varphi : \mathbb{R}^{\Nsig\NG} \times \mathbb{R}^{\Nzeta\NG} \rightarrow \mathbb{R}^{\Ny\NG}$.
With the restriction \eqref{eq:YieldFunc}, the subdifferential of $\phic$ that appears in equations \subeq{eq:GoverningEqGeneral}{{2,3}}
can be written as (Figure \ref{fig:IndFunc:subdifferential})
\begin{equation}
    \partial\phic(\sigma,\zeta) = \{ \lambdabar^\top\nabla\varphi(\sigma,\zeta) | \varphi(\sigma,\zeta) \le 0,\;\lambdabar \ge 0,\; \lambdabar^\top\varphi(\sigma,\zeta) = 0\}
    \label{eq:phic-subgradient}
\end{equation}
\nomenclature[$C$]{$\mathcal{C}$}{Convex elastic region}
\nomenclature[$~$]{$\sqcup_\mathcal{C}$}{Indicator function of the convex elastic region $\mathcal{C}$}
\nomenclature[$\phi~$]{$\varphi$}{Yield function}
\nomenclature{$\lambdabar$, $\lambdabarnp$}{Vector in equations \eqref{eq:phic-subgradient} and \eqref{eq:GoverningEqLagMult} $(\in \mathbb{R}^{\Ny})$;
                                            corresponds to equivalent plastic strain rate; superscript denotes time increment index}
\nomenclature{$x$}{Generic variable used in Figure \ref{fig:IndFunc}}

For quasi-static analyses, the first two terms of equation \subeq{eq:GoverningEqGeneral}{1} are absent.
Also using the fact that $\psic$ has been restricted to being smooth and equation \eqref{eq:phic-subgradient},
equations \eqref{eq:GoverningEqGeneral} can be written in the alternate form
\begin{equation}
  \begin{gathered}
    \begin{alignedat}{4}
        &&
        &&
        &{}\phantom{+}{}B^\top\sigma&
        &=p \\
        &\od{}{t}\nabla_1\psic(\sigma,\zeta)&
        &{}+{}\lambdabar^\top\nabla_1\varphi(\sigma,\zeta)&
        &{}-{}Bv&
        &= 0 \\
        &\od{}{t}\nabla_2\psic(\sigma,\zeta)&
        &{}+{}\lambdabar^\top\nabla_2\varphi(\sigma,\zeta)&
        &&
        &= 0
    \end{alignedat} \\
    \varphi(\sigma,\zeta) \ge 0,\quad \lambdabar \ge 0,\quad \lambdabar^\top\varphi(\sigma,\zeta) = 0
  \end{gathered}
  \label{eq:GoverningEqLagMult}
\end{equation}

Equations \eqref{eq:GoverningEqLagMult} may be formally discretized in time.
Using Backward Euler with a time increment $\Dt$, we get
\begin{equation}
  \begin{gathered}
    \begin{alignedat}{4}
        &&
        &&
        &{}\phantom{+}{}B^\top\sigmanp&
        &=\pnp \\
        &\frac{\nabla_1\psic(\sigmanp,\zetanp)-\nabla_1\psic(\sigman,\zetan)}{\Dt}&
        &{}+{}(\lambdabarnp)^\top\nabla_1\varphi(\sigmanp,\zetanp)&
        &{}-{}B\frac{\unp-\un}{\Dt}&
        &= 0 \\
        &\frac{\nabla_2\psic(\sigmanp,\zetanp)-\nabla_2\psic(\sigman,\zetan)}{\Dt}&
        &{}+{}(\lambdabarnp)^\top\nabla_2\varphi(\sigmanp,\zetanp)&
        &&
        &= 0
    \end{alignedat} \\
    \varphi(\sigmanp,\zetanp) \ge 0,\quad \lambdabarnp \ge 0,\quad (\lambdabarnp)^\top\varphi(\sigmanp,\zetanp) = 0
  \end{gathered}
  \label{eq:GoverningEq-discrete}
\end{equation}
where $u \in \mathbb{R}^{\NDOF}$ is the vector of free displacements.
Multiplying equations \subeq{eq:GoverningEq-discrete}{{2,3}} by $\Dt$, rearranging,
and introducing the symbols $\mu \colonequals \unp-\un$ and $\lambda \colonequals \Dt\lambdabar$,
equations \eqref{eq:GoverningEq-discrete} become
\begin{equation}
  \begin{gathered}
    \begin{alignedat}{4}
        &\nabla_1\psic(\sigmanp,\zetanp)&
        &{}+{}\lambda^\top\nabla_1\varphi(\sigmanp,\zetanp)&
        &{}-{}(B\mu + \underbrace{\nabla_1\psic(\sigman,\zetan)}_{\equalscolon\bsigma})&
        &= 0 \\
        &\nabla_2\psic(\sigmanp,\zetanp)&
        &{}+{}\lambda^\top\nabla_2\varphi(\sigmanp,\zetanp)&
        &{}-{}\underbrace{\nabla_2\psic(\sigman,\zetan)}_{\equalscolon\bzeta}&
        &= 0
    \end{alignedat} \\
    B^\top\sigmanp = \pnp \\
    \varphi(\sigmanp,\zetanp) \ge 0,\quad \lambda \ge 0,\quad \lambda^\top\varphi(\sigmanp,\zetanp) = 0
  \end{gathered}
  \label{eq:GoverningEq-discrete-KKT}
\end{equation}
Equations \eqref{eq:GoverningEq-discrete-KKT} become the starting point for optimization problems of different forms as discussed in the next section.
We note that prescribed displacements can be accounted for by simply replacing all occurrences of the term $(B \mu + \bsigma)$
above by $(B \mu + \Bprsc \muprsc + \bsigma)$, where $\muprsc$ is the increment of prescribed displacements,
and $\Bprsc$ is the deformation-displacement displacement matrix associated with DOF with prescribed displacement.
\nomenclature{$\muprsc$}{Increment of prescribed displacement}
\nomenclature[$Bp$]{$\Bprsc$}{Deformation-displacement matrix associated with DOF with prescribed displacement}

\section{Optimization problems}\label{sec:Opt}
\subsection{Primal optimization problem}\label{sec:Opt:primal}
We recognize immediately that equations \eqref{eq:GoverningEq-discrete-KKT} are the
Karush-Kuhn-Tucker optimality (KKT) conditions \cite{NocedalWright2006} for the following minimization problem.
\begin{equation}
    \begin{gathered}
        (\sigma^{n+1},\zeta^{n+1}) = \argmin_{(\sigma,\zeta)}\ \psic(\sigma,\zeta) - \bsigma^\top\sigma - \bzeta^\top\zeta \\
        \begin{aligned}
            \sbt & B^\top\sigma = \pnp \\
                 & \varphi(\sigma,\zeta) \le 0
        \end{aligned}
    \end{gathered}
    \label{eq:GoverningEq-discrete-opt}
\end{equation}
where $\bsigma$ and $\bzeta$ are as defined in equation \eqref{eq:GoverningEq-discrete-KKT}.
The incremental displacements are the Lagrange multipliers corresponding to the equilibrium constraints (motivating the notation $\mu$).
The Lagrange multipliers $\lambda$ corresponding to the yield constraints can be interpreted as incremental equivalent plastic strains.

\eqref{eq:GoverningEq-discrete-opt} is a convex minimization problem,
since the objective function and inequality constraints are convex, and the equality constraint is linear \cite{Boyd2004}.
We refer to this as the \emph{primal} optimization problem.
We note that \eqref{eq:GoverningEq-discrete-opt} corresponds to the mathematical program referred to as ``approach 1'' in \cite{SivaselvanNLM2013}.
It is identical to equation (46) of reference \cite{Krabbenhoft-InteriorPoint-2007} (where the notation $\kappa$ is used for $\zeta$), and similar to equation (5) in \cite{Bilotta2012}.
It is also the same as equation (17) in \cite{Sivaselvan2009}, where dynamics is considered as well (and the notation $(F,\zeta)$ is used for $(\sigma,\zeta)$).
We assume in the following that
\begin{enumerate}
    \item \eqref{eq:GoverningEq-discrete-opt} has a solution.
    \item $\psic$ is strongly convex, so that the solution of \eqref{eq:GoverningEq-discrete-opt} is unique, and the $\min$ in equation \eqref{eq:DualFn} is meaningful.
    \item \eqref{eq:GoverningEq-discrete-opt} satisfies Slater's constraint qualification condition \cite{Boyd2004}, i.e. that it has a strictly feasible point.
\end{enumerate}
\nomenclature{$n$}{Time increment index}
\nomenclature[$\Deltat$]{$\Dt$}{Time increment}
\nomenclature[$bsigma$]{$\bsigma$, $\bsigmam$}{Defined in equation \eqref{eq:GoverningEq-discrete-KKT}; components of $\bsigma$ corresponding to integration point $m$}
\nomenclature[$bzeta$]{$\bzeta$, $\bzetam$}{Defined in equation \eqref{eq:GoverningEq-discrete-KKT}; components of $\bzeta$ corresponding to integration point $m$}
\nomenclature{$\mu$}{Incremental displacement, $\unp-\un$}
\nomenclature{$\lambda$, $\lambda_m$}{Vector of equivalent plastic strain increments for the entire model $(\in \mathbb{R}^{\Ny\NG})$,
                                      or for a single material point $(\in \mathbb{R}^{\Ny})$}

We write \eqref{eq:GoverningEq-discrete-opt} in a more convenient form
by noting that the objective function and the constraints are \emph{separable} over integrations points.
The objective function can be written as
\begin{equation}
    \Pic(\sigma,\zeta) \colonequals \sumM \psic(\sigmam,\zetam) - \bsigmam^\top\sigmam - \bzetam^\top\zetam
    \label{eq:Pic}
\end{equation}
and the optimization problem as
\begin{equation}
    \begin{gathered}
        (\sigma^{n+1},\zeta^{n+1}) = \argmin_{(\sigma,\zeta)}\ \Pic(\sigma,\zeta) \\
        \begin{aligned}
            \sbt & \sumM B_m^\top\sigmam = p \\
                 & \varphi(\sigmam,\zetam) \le 0,\; \forM
        \end{aligned}
    \end{gathered}
    \label{eq:Primal}
    \tag{PRIMAL}
\end{equation}
where subscript $m$ denotes components of a vector or rows of a matrix corresponding to integration point $m$.
We have used $\psic$ and $\varphi$ in the above equations instead of $\psic_m$ and $\varphi_m$, to minimize notational clutter.
Models where the material is not homogenous, so that these functions are different at different integration points
presents no additional difficulty. We have also dropped the superscript $n+1$ on the load vector $p$ for brevity.
\nomenclature{$m$}{Integration point index $\in \{1,\ldots,\NG\}$}
\nomenclature{$B_m$}{Rows of the matrix $B$ corresponding to integration point $m$ ($\in \mathbb{R}^{\Nsig\times\NDOF}$)}
\nomenclature{$\Pic$}{Objective function of the primal optimization problem (equation \eqref{eq:Pic})}
\nomenclature[$L$]{$\calL$}{Lagrangian of the primal optimization problem (equation \eqref{eq:Lagrangian})}

The Lagrangian associated with \eqref{eq:Primal} is
\begin{equation}
    \begin{aligned}
        \calL\del[1]{\sigzetavec,\mulamvec} &= \Pic(\sigma,\zeta) - \mu^\top\left(\sumM B^\top_m \sigma_m - p\right) + \sumM\lambdam^\top \varphi(\sigmam, \zetam) \\
                                          &= \mu^\top p + \sumM\sbr[1]{\psic(\sigmam, \zetam) + \lambdam^\top\varphi(\sigmam,\zetam) - (B_m \mu + \bsigmam)^\top\sigmam - \bzetam^\top\zetam} \\
                                          &= \mu^\top p + \sumM \calL_m{\sigzetavecm,\mulamvecm}
    \end{aligned}
    \label{eq:Lagrangian}
\end{equation}
where we have defined
\begin{equation}
    \calL_m\del[1]{\sigzetavecm,\mulamvecm} \colonequals \psic(\sigmam, \zetam) + \lambdam^\top\varphi(\sigmam,\zetam) - (B_m \mu + \bsigmam)^\top\sigmam - \bzetam^\top\zetam
    \label{eq:LagrangianIPm}
\end{equation}
The Karush-Kuhn-Tucker (KKT) optimality conditions for \eqref{eq:Primal} follow from the Lagrangian \eqref{eq:Lagrangian} as
\begin{equation}
    \begin{gathered}
        \begin{alignedat}{6}
            \LTosigmam = 0&:&&\Dif_1\psic(\sigmam,\zetam)&&{}+{}\lambdam^\top\Dif_1\varphi(\sigmam,\zetam)&&{}-{}(B_m\mu+\bsigmam)^\top&&=0&&,\;\forM \\
            \LTozetam = 0   &:&&\Dif_2\psic(\sigmam,\zetam)&&{}+{}\lambdam^\top\Dif_2\varphi(\sigmam,\zetam)&&{}-{}\bzetam^\top&            &=0&&,\;\forM
        \end{alignedat} \\
        \sumM B_m^\top \sigmam = p \\
        \varphi(\sigmam,\zetam) \leq 0,\;\lambdam \ge 0,\;\lambdam^\top\varphi(\sigmam,\zetam) = 0,\;\forM
    \end{gathered}
    \label{eq:Primal-KKT}
    \tag{KKT}
\end{equation}
These are merely equations \eqref{eq:GoverningEq-discrete-KKT} separated over the integration points.
Since the gradient is the transpose of the derivative, while equations \subeq{eq:GoverningEq-discrete-KKT}{{1,2}} are in terms of column vectors,
equations \subeq{eq:Primal-KKT}{{1,2}} are in terms of row vectors.
We next construct the dual of problem \eqref{eq:Primal}, which in section \ref{sec:reformulate} leads to a reformulation of the classical displacement-based approach as an optimization problem.
\nomenclature[$Lm$]{$\calL_m$}{Contribution to primal Lagrangian of integration point $m$ (equation \eqref{eq:LagrangianIPm})}

\subsection{Dual optimization problem}\label{sec:Opt:dual}
Duality is familiar in mechanics from the classical principles of total potential energy and total complementary potential energy,
and from the upper- and lower-bound theorems of limit analysis.
Associated with every optimization problem is a dual problem, which often has a useful interpretation.
For convex problems such as \eqref{eq:Primal} that satisfy constraint qualification conditions as assumed above, \emph{strong duality} holds.
This means that a dual optimal solution is equal in objective function value to a primal optimal solution, and that the primal-dual solution pair satisfies the KKT conditions.
The KKT conditions are also sufficient conditions for a solution of the primal problem \cite{Boyd2004}.

To construct the dual problem for \eqref{eq:Primal}, we start from the Lagrange dual function \cite{Boyd2004},
defined as the minimum of the Lagrangian of $\Pic$ over $(\sigma,\zeta)$.
\begin{equation}
    -\Pi(\mu,\lambda) \colonequals \min_{(\sigma,\zeta)} \calL\del[1]{\sigzetavec,\mulamvec}
    \label{eq:DualFn}
\end{equation}
We are able to write $\min$ here instead of $\inf$, because of the strong convexity assumption on $\psic$.
Let $(\sigmstar,\zetamstar), \forM$ be the minimizers in \eqref{eq:DualFn}.
Functions $\sigmstar(\mu,\lambda)$ and $\zetamstar(\mu,\lambda)$ are then implicitly defined by the
following optimality conditions.
\begin{equation}
    \begin{alignedat}{4}
        &\Dif_1\psic(\sigmstar(\mu,\lambda),\zetamstar(\mu,\lambda))&&{}+{}\lambdam^\top\Dif_1\varphi(\sigmstar(\mu,\lambda),\zetamstar(\mu,\lambda))&&{}-{}(B_m\mu+\bsigmam)^\top&&=0 \\
        &\Dif_2\psic(\sigmstar(\mu,\lambda),\zetamstar(\mu,\lambda))&&{}+{}\lambdam^\top\Dif_2\varphi(\sigmstar(\mu,\lambda),\zetamstar(\mu,\lambda))&&{}-{}\bzetam^\top&            &=0
    \end{alignedat} \\
    \label{eq:Def-SigZetaStar}
\end{equation}
i.e., given $(\mu,\lambda)$, the functions $\sigmstar(\mu,\lambda)$ and $\zetamstar(\mu,\lambda)$ can be evaluated by solving equation \eqref{eq:Def-SigZetaStar}.
Clearly, this is equivalent to unconstrained minimization of the Lagrangian function of equation \eqref{eq:LagrangianIPm}.
\nomenclature{$\Pi$}{Objective function of dual problem, defined in equation \eqref{eq:DualFn}, explicit formula in equation \eqref{eq:DualFn-explicit}}
\nomenclature{$\sigmstar$}{Minimizer of equation \eqref{eq:DualFn} defined implicity in equations \eqref{eq:Def-SigZetaStar}}
\nomenclature{$\zetamstar$}{Minimizer of equation \eqref{eq:DualFn} defined implicity in equations \eqref{eq:Def-SigZetaStar}}

The Lagrange dual function $\Pi$ of equation \eqref{eq:DualFn} may be written explicitly as
\begin{equation}
    \begin{aligned}
        \Pi(\mu,\lambda) = -\mu^\top p - \sumM &\left[\psic(\sigmstar(\mu,\lambda), \zetamstar(\mu,\lambda))
                                                       + \lambdam^\top\varphi(\sigmstar(\mu,\lambda),\zetamstar(\mu,\lambda))
                                                       \vphantom{- (B_m \mu + \bsigmam)^\top\sigmstar(\mu,\lambda)
                                                       - \bzetam^\top\zetamstar(\mu,\lambda)} \right. \\
                                               &\left. - (B_m \mu + \bsigmam)^\top\sigmstar(\mu,\lambda)
                                                       - \bzetam^\top\zetamstar(\mu,\lambda)\right]
    \end{aligned}
    \label{eq:DualFn-explicit}
\end{equation}
The dual optimization problem corresponding to \eqref{eq:Primal} is then
\begin{equation}
    \begin{gathered}
        \min \Pi(\mu,\lambda) \\
        \sbt\ \lambda \ge 0
    \end{gathered}
    \label{eq:Dual}
    \tag{DUAL}
\end{equation}
As discussed above, solving this problem is completely equivalent to solving \eqref{eq:Primal}.
\eqref{eq:Dual} is a convex optimization problem with \emph{simple bound constraints}.
It is closely related to the mathematical program referred to as approach 2 in \cite{SivaselvanNLM2013}.
In section \ref{sec:Opt:reduced} next, we reduce this to an unconstrained optimization problem in $\mu$ alone.
This forms the basis for reformulation of the classical displacement-based approach in section \ref{sec:reformulate}.

\subsection{Reduced dual problem}\label{sec:Opt:reduced}
Due to the separability of the objective function $\Pi$, minimization over $\lambda$ in \eqref{eq:Dual} may be carried into the summation in \eqref{eq:DualFn-explicit}.
Thus defining the function
\begin{equation}
    \begin{aligned}
        \bar{\Pi}(\mu) \colonequals -\mu^\top p
                         - \sumM \min_{\lambdam \ge 0} &\left[\psic(\sigmstar(\mu,\lambda), \zetamstar(\mu,\lambda))
                                                                + \lambdam^\top\varphi(\sigmstar(\mu,\lambda),\zetamstar(\mu,\lambda))
                                                                \vphantom{- (B_m \mu + \bsigmam)^\top\sigmstar(\mu,\lambda)
                                                                - \bzetam^\top\zetamstar(\mu,\lambda)} \right. \\
                                                         &\left. - (B_m \mu + \bsigmam)^\top\sigmstar(\mu,\lambda)
                                                                - \bzetam^\top\zetamstar(\mu,\lambda)\right]
    \end{aligned}
    \label{eq:Pibar}
\end{equation}
\eqref{eq:Dual} can be reduced to minimization of $\mu$ alone as
\begin{equation}
    \min_\mu \bar{\Pi}(\mu)
    \label{eq:DualProblem-Nested}
    \tag{REDUCED}
\end{equation}
This reduced problem leads to the reformulation of classical displacement-based approach as an optimization problem.
We note that in this way, the solution of \eqref{eq:Dual} is carried out in a nested fashion
--- the minimization over $\lambdam \ge 0$ within the summation is at the integration point level,
while the minimization over $\mu$ is at the global level.
We reiterate that the integration point-level minimization format is well-known in the context of so-called return-mapping schemes,
discussed further in section \ref{sec:return-map}.
However, the minimization structure at the global level \eqref{eq:DualProblem-Nested} is not commonly recognized.
We take on the solution of \eqref{eq:DualProblem-Nested} in section \ref{sec:reformulate}.
In the remainder of this section, we develop an explicit expression for $\bar{\Pi}$.

Consider the minimization within the summation in equation \eqref{eq:Pibar}.
\begin{equation}
    \begin{gathered}
        \begin{aligned}
            \min_{\lambdam}\; &\psic(\sigmstar(\mu,\lambda), \zetamstar(\mu,\lambda)) + \lambdam^\top\varphi(\sigmstar(\mu,\lambda),\zetamstar(\mu,\lambda)) \\
                              & - (B_m \mu + \bsigmam)^\top\sigmstar(\mu,\lambda) - \bzetam^\top\zetamstar(\mu,\lambda)
        \end{aligned} \\
        \sbt \lambdam \ge 0
    \end{gathered}
    \label{eq:IP-dual}
    \tag{IPDUAL}
\end{equation}
where the prefix IP is used in the equation tag since the minimization is at the integration point level.
Let $\lamtilm(\mu)$ be the minimizer. We recognize that \eqref{eq:IP-dual} is the Lagrange dual of
\begin{equation}
    \begin{gathered}
        \min_{\sigma,\zeta}\; \psic(\sigma,\zeta) - (B_m\mu + \bsigmam)^\top\sigma - \bzetam^\top\zeta \\
        \sbt \varphi(\sigma,\zeta) \le 0
    \end{gathered}
    \label{eq:IP-primal}
    \tag{IPPRIMAL}
\end{equation}
Let $(\sigtildem(\mu),\zetatildem(\mu))$ be the minimizer, given implicitly by the optimality conditions
\begin{equation}
    \begin{gathered}
        \begin{alignedat}{4}
            &\Dif_1\psic(\sigtildem(\mu),\zetatildem(\mu))&&{}+{}\lamtilm(\mu)^\top\Dif_1\varphi(\sigtildem(\mu),\zetatildem(\mu))&&{}-{}(B_m\mu + \bsigmam)^\top&&= 0 \\
            &\Dif_2\psic(\sigtildem(\mu),\zetatildem(\mu))&&{}+{}\lamtilm(\mu)^\top\Dif_2\varphi(\sigtildem(\mu),\zetatildem(\mu))&&{}-{}\bzetam^\top&& = 0
        \end{alignedat} \\
        \varphi(\sigtildem(\mu),\zetatildem(\mu)) \le 0,\quad \lamtilm(\mu) \ge 0,\quad \lamtilm(\mu)^\top\varphi(\sigtildem(\mu),\zetatildem(\mu)) = 0
    \end{gathered}
    \label{eq:IP-primal-KKT}
\end{equation}
By strong duality, the optimum values of \eqref{eq:IP-dual} and \eqref{eq:IP-primal} are identical.
Therefore, $\bar{\Pi}$ can be written as
\begin{equation}
    \bar{\Pi}(\mu) = -\mu^\top p - \sumM \left[\psic(\sigtildem(\mu),\zetatildem(\mu)) - (B_m\mu + \bsigmam)^\top\sigtildem(\mu) - \bzetam^\top\zetatildem(\mu)\right]
    \label{eq:Pibar-alt}
\end{equation}
We note as an aside that $(\sigtildem(\mu),\zetatildem(\mu),\lamtilm(\mu))$ solve the KKT conditions of \eqref{eq:IP-primal},
and that the identities $\sigmstar(\mu,\lamtilm(\mu)) = \sigtildem(\mu)$ and $\zetamstar(\mu,\lamtilm(\mu)) = \zetatildem(\mu)$ hold.
Equation \eqref{eq:Pibar-alt} provides an explicit expression for $\bar{\Pi}$.
In the next section, we present an algorithm for its minimization, i.e. to solve \eqref{eq:DualProblem-Nested}.
\nomenclature[$\Pibar$]{$\bar{\Pi}$}{Objective function of the reduced dual problem defined in equation \eqref{eq:Pibar}}
\nomenclature[$\lambdatilde$]{$\lamtilm$}{Minimizer of \eqref{eq:IP-dual}}
\nomenclature[$\sigmatildem$]{$\sigtildem$}{Minimizer of \eqref{eq:IP-primal}}
\nomenclature[$\zetatildem$]{$\zetatildem$}{Minimizer of \eqref{eq:IP-primal}}

\section{Optimization reformulation of the classical displacement-based nested approach}\label{sec:reformulate}
In this section, we consider numerical solution of the reduced dual problem \eqref{eq:DualProblem-Nested}.
We seek to use Newton's method for unconstrained optimization.
For this, we require the first and second derivatives of the objective function \eqref{eq:Pibar-alt}.
We discuss differentiability of this function, and computation of these derivatives in section \ref{sec:reformulate:derivatives} below.
We then present Newton's method in section \ref{sec:reformulate:algorithm}.
Solving the reduced dual problem in this manner may be interpreted as a reformulation of the classical displacement-based nested approach as an optimization problem.
The reason for this interpretation is also discussed in section \ref{sec:reformulate:derivatives}.

\subsection{Differentiability and derivatives}\label{sec:reformulate:derivatives}
We seek to obtain the first and second derivatives of $\bar{\Pi}$ needed for minimization using Newton's method.
In computing these derivatives, is helpful to think of \eqref{eq:IP-primal} as a parametric convex optimization problem with parameter $\mu$.
The first derivative of $\bar{\Pi}$ is well-defined, and can be obtained by differentiating \eqref{eq:Pibar-alt} as
\begin{equation}
    \nabla\bar{\Pi}(\mu) = \sumM B_m^\top\sigtildem(\mu) - p
    \label{eq:Pibar-first-deriv}
\end{equation}
The optimality condition for \eqref{eq:DualProblem-Nested},
$\nabla\bar{\Pi}(\mu) = 0$, is thus
\begin{equation}
    \sumM B_m^\top\sigtildem(\mu) = p
    \label{eq:equilibrium}
\end{equation}
which is the equilibrium equation.

Equation \eqref{eq:equilibrium} can also be obtained directly from \eqref{eq:GoverningEq-discrete-KKT},
since \subeq{eq:GoverningEq-discrete-KKT}{1,2,4} are the KKT conditions of \eqref{eq:IP-primal}.
When derived this way however, the minimization structure \eqref{eq:DualProblem-Nested} at the global level is not obvious.
Recognizing this structure enables direct application of Newton's method for unconstrained minimization in section \ref{sec:reformulate:algorithm},
without need for heuristics, often used for example in step-length determination.
The classical displacement-based nested approach used \eqref{eq:equilibrium} as the starting point, rather than the minimization problem \eqref{eq:DualProblem-Nested}.
It is in this sense that we consider \eqref{eq:DualProblem-Nested} as a reformulation of the classical approach.

To compute the second derivative of $\bar{\Pi}$, we need the derivative of the functions $\sigtildem(\mu)$.
However, in trying to differentiate \eqref{eq:IP-primal-KKT} to obtain this, we realize that
$\sigtildem(\mu)$ and $\zetatildem(\mu)$ are not differentiable when they are degenerate solutions of \eqref{eq:IP-primal},
i.e., when $\varphi(\sigtildem(\mu),\zetatildem(\mu)) = 0$ and $\lamtilm(\mu) = 0$\footnote{Indeed
these degenerate points feature in the context of the first derivative \eqref{eq:Pibar-first-deriv} as well,
where an argument has to be made that $\lamtilm(\mu)^\top \Dif\varphi = 0$ at such points.}.
In fact at such points, these functions are only directionally differentiable \cite{DempeSolutionDifferentiability}.
When the minimizer of \eqref{eq:DualProblem-Nested} is such a degenerate point,
the rate of convergence of Newton's method is only super-linear rather than quadratic \cite{Qi1993}.

Here, we proceed formally ignoring these degenerate points, the set of which has measure zero in $\mathbb{R}^{\NDOF}$.
We first introduce the index set
\begin{equation}
    \calI_m = \left\{k \in \{1,\ldots,\Ny\} \; \big| \; (\lambda_m(\mu))_k > 0 \, \right\}
    \label{eq:index-set}
\end{equation}
The optimality conditions \eqref{eq:IP-primal-KKT} can then be rewritten as
\begin{equation*}
    \begin{gathered}
        \begin{alignedat}{4}
            &\Dif_1\psic(\sigtildem(\mu),\zetatildem(\mu))&&{}+{}\sum_{k\in\calI_m}(\lamtilm(\mu))_k\Dif_1\varphi_k(\sigtildem(\mu),\zetatildem(\mu))&&{}-{}(B_m\mu + \bsigmam)^\top&&= 0 \\
            &\Dif_2\psic(\sigtildem(\mu),\zetatildem(\mu))&&{}+{}\sum_{k\in\calI_m}(\lamtilm(\mu))_k\Dif_2\varphi_k(\sigtildem(\mu),\zetatildem(\mu))&&{}-{}\bzetam^\top&& = 0
        \end{alignedat} \\
        \varphi_k(\sigtildem(\mu),\zetatildem(\mu)) = 0,\; \text{for}\ k \in \calI_m
    \end{gathered}
\end{equation*}
Differentiating this, we get
\begin{subequations}
\label{eq:IP-derivative}
\begin{equation}
    \left[ \begin{array}{c|c|c}
        \Kss & \Ksx & \Phis^\top \\
        \hline
        \Ksx^\top & \Kxx & \Phix^\top \\
        \hline
        \Phis & \Phix & 0
    \end{array} \right]
    \left( \begin{array}{c}
        \Dif\sigtildem(\mu) \\ \hline \Dif\zetatildem(\mu) \\ \hline \Dif(\lamtilm(\mu)_{\calI_m :})
    \end{array} \right)
    =
    \left( \begin{array}{c}
        B_m \\ \hline 0 \\ \hline 0
    \end{array} \right)
    \label{eq:IP-derivative:lineq}
\end{equation}
where
\begin{equation}
    \begin{aligned}
        \Kss &= \Dif_1^2\psic(\sigtildem(\mu),\zetatildem(\mu)) + \sum_{k\in\calI_m}(\lamtilm(\mu))_k\Dif_1^2\varphi_k(\sigtildem(\mu),\zetatildem(\mu)) \\
        \Ksx &= \Dif_2\Dif_1\psic(\sigtildem(\mu),\zetatildem(\mu)) + \sum_{k\in\calI_m}(\lamtilm(\mu))_k\Dif_2\Dif_1\varphi_k(\sigtildem(\mu),\zetatildem(\mu)) \\
        \Kxx &= \Dif_2^2\psic(\sigtildem(\mu),\zetatildem(\mu)) + \sum_{k\in\calI_m}(\lamtilm(\mu))_k\Dif_2^2\varphi_k(\sigtildem(\mu),\zetatildem(\mu)) \\
        \Phis &= \del{\Dif_1\varphi(\sigtildem(\mu),\zetatildem(\mu))}_{\calI_m :} \\
        \Phix &= \del{\Dif_2\varphi(\sigtildem(\mu),\zetatildem(\mu))}_{\calI_m :}
    \end{aligned}
    \label{eq:IP-derivative:defs}
\end{equation}
In the last two equations above, we have used MATLAB-style indexing \cite{Mathworks}.
The subscript $\calI_m:$ extracts the submatrix consisting of rows $\calI_m$ and all the columns.
\end{subequations}
\nomenclature[$I_m$]{$\calI_m$}{Index set of activated yield conditions defined in equation \eqref{eq:index-set}}
\nomenclature{$k$}{Yield function component index within an integration point $\in \{1,\ldots,\Ny\}$}
\nomenclature[$K_\sigma$]{$\Kss$, $\Ksx$, $\Kxx$}{Derivative matrices defined in equation \eqref{eq:IP-derivative:defs}}
\nomenclature[$\Phi_\sigma$]{$\Phis$, $\Phix$}{Derivative matrices defined in equation \eqref{eq:IP-derivative:defs}}

\begin{subequations}
The derivative $\Dif\sigtildem(\mu)$ can be obtained by solving \eqref{eq:IP-derivative:lineq}.
It is of the form $\Kbar_m B_m$, where $\Kbar_m$ is referred to as the \emph{consistent tangent stiffness} matrix \cite{SimoTaylor1985}.
When the submatrices in \eqref{eq:IP-derivative:lineq} have a simple form, it may also be less expensive to compute the consistent tangent directly as
\label{eq:Kbar-direct}
\begin{equation}
    \Kbar_m = \Kss^{-1}\left[\Kss + \Ksx\Kbarxx^{-1}\Ksx^\top - \Phibars^\top\Xi^{-1}\Phibars\right]\Kss^{-1}
    \label{eq:Kbar-direct:Kbar}
\end{equation}
where
\begin{equation}
    \begin{aligned}
        \Kbarxx  &= \Kxx - \Ksx^\top\Kss^{-1}\Ksx \\
        \Phibarx &= \Phix - \Phis\Kss^{-1}\Ksx \\
        \Phibars &= \Phis - \Phibarx\Kbarxx^{-1}\Ksx^\top \\
        \Xi      &= \Phis\Kss^{-1}\Phis^\top + \Phibarx\Kbarxx^{-1}\Phibarx^\top
    \end{aligned}
    \label{eq:Kbar-direct:defs}
\end{equation}
\end{subequations}
\nomenclature[$Kbar$]{$\Kbar_m$}{Consistent tangent stiffness matrix for integration point $m$}
\nomenclature[$Kbar_\zeta$]{$\Kbarxx$}{Derivative matrix defined in equation \eqref{eq:Kbar-direct:defs}}
\nomenclature[$Phibar_\sigma$]{$\Phibars$, $\Phibarx$}{Derivative matrices defined in equation \eqref{eq:Kbar-direct:defs}}
\nomenclature{$\Xi$}{Derivative matrix defined in equation \eqref{eq:Kbar-direct:defs}}

Now, by differentiating equation \eqref{eq:Pibar-first-deriv}, the second derivative of the objective function $\bar{\Pi}$,
the tangent stiffness matrix, may be written as
\begin{equation}
    \nabla^2\bar{\Pi}(\mu) = \sumM B_m^\top \Kbar_m B_m
    \label{eq:tangent}
\end{equation}

\subsection{Optimization algorithm}\label{sec:reformulate:algorithm}
With the first derivative $\nabla\bar{\Pi}(\mu)$, the vector of unbalanced forces,
and second derivative $\nabla^2\bar{\Pi}(\mu)$ the tangent stiffness matrix, at hand,
the energy $\bar{\Pi}(\mu)$ can be minimized
using the standard Newton's method for unconstrained minimization \cite{Boyd2004,NocedalWright2006},
outlined in Procedure \ref{alg:Newt}.
A Newton direction $\dmu$ is computed by solving the linear system
\begin{equation}
    \nabla^2\bar{\Pi}(\mu^\iter)\dmu = -\nabla\bar{\Pi}(\mu^\iter)
    \label{eq:search-direction}
\end{equation}
To assure convergence for large increments, a line search is typically needed.
A step length $s$ is computed by backtracking to satisfy the following sufficient decrease condition.
\begin{equation}
    \bar{\Pi}(\mu + s\dmu) \le \bar{\Pi}(\mu) + \beta s \nabla\bar{\Pi}(\mu)^\top\dmu
    \label{eq:suff-decr}
\end{equation}
where $\beta$ is a sufficient decrease parameter to be in the interval $(0,0.5)$ \cite{NocedalWright2006}, taken in the example in section \ref{sec:NumExample} as $10^{-4}$.
Solving the problem as one of optimization removes the need to use any other heuristics to ensure convergence.
The working of the backtracking procedure can be visualized as shown in Figure \ref{fig:linesearch} in section \ref{sec:NumExample}.
\begin{algorithm}
    \begin{algorithmic}[1]
        \State Given force $p$, prescribed displacement $\muprsc$, $(\bsigmam,\bzetam)$ for each integration point $m$, and starting point $\mu^0$
        \For{$\iter \gets 0, \maxiter$}
            \State Compute $\bar{\Pi}(\mu^\iter)$, $\nabla\bar{\Pi}(\mu^\iter)$ and $\nabla^2\bar{\Pi}(\mu^\iter)$ using Procedure \ref{alg:func-grad-hess}
            \If{$||\nabla\bar{\Pi}(\mu^\iter)|| \le \mathtt{tol}$}
                \State \textbf{break}
            \EndIf
            \State Compute search direction $\dmu$ by solving \eqref{eq:search-direction}
            \State $s = 1$
            \For{$\nback \gets 1, \maxbacktrack$} \Comment{backtracking line search}
                \State Compute $\bar{\Pi}(\mu^\iter + s\dmu)$ using Procedure \ref{alg:func-grad-hess}
                \If{sufficient decrease \eqref{eq:suff-decr} is obtained}
                    \State \textbf{break}
                \EndIf
                \State $s \gets \gamma s$ \Comment{$\gamma$ is a backtracking parameter}
            \EndFor
            \State $\mu^{\iter+1} \gets \mu^\iter + s\dmu$
        \EndFor
    \end{algorithmic}
    \caption{\textsc{Newton's method for $\min \bar{\Pi}(\mu)$}}
    \label{alg:Newt}
\end{algorithm}
\nomenclature[$\delta\mu$]{$\dmu$}{Newton search direction (equation \eqref{eq:search-direction})}
\nomenclature[$MAXITER$]{$\maxiter$}{Maximum number of Newton iterations in Procedure \ref{alg:Newt}}
\nomenclature[$tol$]{$\mathtt{tol}$}{Convergence tolerances for Newton's method in Procedure \ref{alg:Newt}}
\nomenclature[$iter$]{$\iter$}{Iteration count in Newton's method in Procedure \ref{alg:Newt}}
\nomenclature[$\muiter$]{$\mu^\iter$}{$\mu$ in iteration number $\iter$ of Newton's method (Procedure \ref{alg:Newt})}
\nomenclature{$s$}{Length of Newton step in \ref{alg:Newt}}
\nomenclature[$nback$]{$\nback$}{Number of backtracking steps in Procedure \ref{alg:Newt}}
\nomenclature[$MAXBACKTRACK$]{$\maxbacktrack$}{Maximum number of backtracking steps in Procedure \ref{alg:Newt}}
\nomenclature{$\beta$}{Sufficient decrease parameter used in choice of step length (equation \eqref{eq:suff-decr})}
\nomenclature{$\gamma$}{Backtracking parameter used in Procedure \ref{alg:Newt}}

\begin{algorithm}
    \begin{algorithmic}[1]
        \State Given $\mu$, $p$, $\muprsc$
        \State $\bar{\Pi} \gets -\mu^\top p$ ; $\nabla\bar{\Pi} \gets -p$ ; $\nabla^2\bar{\Pi} \gets 0$
        \For{$m \gets 1, \NG$}
            \State\label{alg:func-grad-hess:IP-Primal} Compute $(\sigtildem,\zetatildem)$ by solving \eqref{eq:IP-primal}
            \State $\bar{\Pi} \gets \bar{\Pi} - \left[\psic(\sigtildem,\zetatildem) - (B_m\mu + \bsigmam)^\top\sigtildem - \bzetam^\top\zetatildem\right]$
            \State Assemble $B_m^\top \sigtildem$ into $\nabla\bar{\Pi}$
            \State Computer $\Kbar_m$ using either of equations \eqref{eq:IP-derivative:lineq} or \eqref{eq:Kbar-direct:Kbar}
            \State Assemble $B_m^\top \Kbar_m B_m$ into $\nabla^2\bar{\Pi}$
        \EndFor
    \end{algorithmic}
    \caption{\textsc{Evaluate $\bar{\Pi}(\mu)$, $\nabla\bar{\Pi}(\mu)$ and $\nabla^2\bar{\Pi}(\mu)$}}
    \label{alg:func-grad-hess}
\end{algorithm}

\section{Solution of \eqref{eq:IP-primal} for specific constitutive models}\label{sec:return-map}
In the framework described here, a constitutive model is specified by a complementary stored energy function $\psic$,
and a yield function $\varphi$ or complementary dissipation function $\phic$.
In general, \eqref{eq:IP-primal} may be solved (step \ref{alg:func-grad-hess:IP-Primal} of Procedure \ref{alg:func-grad-hess})
using any standard algorithm for inequality constrained convex optimization such as an interior-point method or a sequential quadratic programming method \cite{Boyd2004,NocedalWright2006}.
Alternately, it may also be advantageous on occasion to solve \eqref{eq:IP-dual}, which has simple bound constraints.
Such a general approach may often be necessary, particularly for more complex models such as multi-surface plasticity models.
However in some cases, specific forms of the complementary stored energy function $\psic$ and the yield function $\varphi$
may facilitate simpler strategies to solve \eqref{eq:IP-primal}. These strategies are known as return-mapping algorithms.
In the following, we discuss one such case, namely an elastoplastic model characterized by isotropic linear elasticity, linear kinematic hardening, nonlinear isotropic hardening,
and von Mises yield condition. This model is described by
\begin{equation}
    \begin{aligned}
        \psic(\sigma,\zeta)   &= \half \sigma^\top C^{-1} \sigma + \half \zetakh^\top H^{-1} \zetakh + \psicih(\zetaih) \\
        \varphi(\sigma,\zeta) &= \sqrt{ (\sigma-\zetakh)^\top P (\sigma-\zetakh) } - \sqrttt(\sigmay + \zetaih)
    \end{aligned}
    \label{eq:elastoplastic}
\end{equation}
Here, $\zetakh$ and $\zetaih$ are components of the internal variable $\zeta$ corresponding to kinematic hardening and isotropic hardening respectively.
$\zetakh$ is often referred to as the back stress, and is the same type of object as $\sigma$. $\zetakh$ is a scalar.
$C$ and $H$ are matrices of elastic and kinematic hardening moduli. The quadratic terms in $\psic$ correspond to linear elasticity and linear kinematic hardening.
The function $\psicih$ is not necessarily quadratic, and represents nonlinear isotropic hardening.
The relationship between this form of $\psic$ and the hardening model discussed in \cite{SimoHughesCompInelasticity} is discussed in \cite{SivaselvanNLM2013}.
von Mises yield function in both 3D and 2D plane-stress cases can be expressed in the form $\varphi$ above with appropriate choice of the matrix $P$ (Table \ref{tab:PQ}) \cite{SimoHughesCompInelasticity}.
$\sigmay$ is the uniaxial yield stress.
\nomenclature{$\zetakh$, $\zetaih$}{Components of internal variable $\zeta$ corresponding to kinematic and isotropic hardening, used in constitutive model equation \eqref{eq:elastoplastic}}
\nomenclature{$C$}{Matrix of elastic moduli in the constitutive model equation \eqref{eq:elastoplastic}}
\nomenclature{$H$}{Matrix of kinematic hardening moduli in the constitutive model equation \eqref{eq:elastoplastic}}
\nomenclature{$\psicih$}{Part of complementary stored energy representing kinematic hardening in the constitutive model equation \eqref{eq:elastoplastic}}
\nomenclature{$P$}{Matrix in von Mises yield condition \subeq{eq:elastoplastic}{2}}
\nomenclature{$\sigmay$}{Uniaxial yield stress in von Mises yield condition \subeq{eq:elastoplastic}{2}}

The optimality conditions \eqref{eq:IP-primal-KKT} corresponding to the model \eqref{eq:elastoplastic} are
\begin{equation}
    \begin{gathered}
        \begin{alignedat}{4}
            &C^{-1}\sigma&     &{}+{}\lambda\frac{P(\sigma-\zetakh)}{\sqrt{ (\sigma-\zetakh)^\top P (\sigma-\zetakh) }}&&{}-{}(B\mu + \bsigma)&&=0 \\
            &H^{-1}\zetakh&    &{}-{}\lambda\frac{P(\sigma-\zetakh)}{\sqrt{ (\sigma-\zetakh)^\top P (\sigma-\zetakh) }}&&{}-{}\bzetakh&        &=0 \\
            &{\psicih}'(\zetaih)&&{}-{}\lambda\sqrttt&                                                                    &{}-{}\bzetaih&        &=0
        \end{alignedat} \\
        \varphi(\sigma,\zeta) \le 0,\quad \lambda \ge 0,\quad \lambda\varphi(\sigma,\zeta) = 0
    \end{gathered}
    \label{eq:elastoplastic-KKT}
\end{equation}
where the prime on $\psicih$ denotes the first derivative.
For brevity, we have written $\sigma$ for $\sigtildem(\mu)$ etc.\ and dropped the subscript $m$.
Since there is a single yield condition, there are only two possibilities.
\begin{enumerate}
    \item $\varphi(\sigma,\zeta) \le 0$. In this case, the material point is elastic in the increment, and $\lambda = 0$ from equation \subeq{eq:elastoplastic-KKT}{4}.
        From equations \subeq{eq:elastoplastic-KKT}{1,2,3}, an \emph{elastic trial state} is computed.
        \begin{equation}
            \sigmatr = \sigman + CB\mu,\quad \zetakhtr = \zetakh^n,\quad \zetaihtr = \zetaih^n
            \label{eq:trial-state}
        \end{equation}
        If $\varphi(\sigmatr,\zetatr) \le 0$ then this is the solution of \eqref{eq:IP-primal}.
    \item If $\varphi(\sigmatr,\zetatr) > 0$, it can be concluded the material point plastifies.
        So $\lambda > 0$, and \subeq{eq:elastoplastic-KKT}{1,2,3,4}
        constitute four equations in the four unknowns $\sigma$, $\zetakh$, $\zetaih$ and $\lambda$.
        To solve these, a further property of the matrices involved may be invoked.
        In both 3D and 2D plane-stress cases, the matrices $C$, $H$ and $P$ share the same eigenvectors.
        Therefore, there is an orthogonal matrix $Q$ (Table \ref{tab:PQ}), consisting of these eigenvectors as columns,
        which simultaneously diagonalizes all three matrices, i.e. $Q^\top P Q$ is diagonal matrix $\Lambda_P$ etc.
        \begin{table}
            \centering
            \begin{tabular}{lcc}
                \toprule
                                    & $P$ & $Q$ \\
                \midrule
                    3D              & $\frac{1}{3}\begin{bmatrix}
                                                     2  & -1 & -1 & 0 & 0 & 0 \\
                                                     -1 &  2 & -1 & 0 & 0 & 0 \\
                                                     -1 & -1 &  2 & 0 & 0 & 0 \\
                                                      0 & 0 & 0 & 6 & 0 & 0 \\
                                                      0 & 0 & 0 & 0 & 6 & 0 \\
                                                      0 & 0 & 0 & 0 & 0 & 6
                                                  \end{bmatrix}$
                                   & $\begin{bmatrix}
                                         \frac{1}{\sqrt{3}} & \sqrt{\frac{2}{3}}  & 0                   & 0 & 0 & 0 \\
                                         \frac{1}{\sqrt{3}} & -\frac{1}{\sqrt{6}} & \frac{1}{\sqrt{2}}  & 0 & 0 & 0 \\
                                         \frac{1}{\sqrt{3}} & -\frac{1}{\sqrt{6}} & -\frac{1}{\sqrt{2}} & 0 & 0 & 0 \\
                                         0 & 0 & 0 & 1 & 0 & 0 \\
                                         0 & 0 & 0 & 0 & 1 & 0 \\
                                         0 & 0 & 0 & 0 & 0 & 1
                                     \end{bmatrix}$ \\
                \midrule
                    2D plane-stress & $\frac{1}{3}\begin{bmatrix} 2 & -1 & 0 \\ -1 & 2 & 0 \\ 0 & 0 & 6 \end{bmatrix}$
                                    & $\begin{bmatrix} \frac{1}{\sqrt{2}} & \frac{1}{\sqrt{2}} & 0 \\ \frac{1}{\sqrt{2}} & -\frac{1}{\sqrt{2}} & 0 \\ 0 & 0 & 1 \end{bmatrix}$ \\
                \bottomrule
            \end{tabular}
            \caption{von Mises yield condition matrix $P$ and diagonalizing matrix $Q$}
            \label{tab:PQ}
        \end{table}
        This has been pointed out for the 2D case in \cite{SimoTaylorPlaneStress}.
        The 3D case is generally presented differently, by decomposing the stress into volumetric and deviatoric parts.
        In both cases, however, $Q$ represents transformation of stress components to volumentric-deviatoric coordinates.
        Equations \subeq{eq:elastoplastic-KKT}{1,2,3,4} can be rewritten as
        \begin{equation}
            \begin{gathered}
                \begin{alignedat}{4}
                    &\sigma& &{}+{}\frac{\lambda}{\sqrttt(\sigmay + \zetaih)}CP(\sigma-\zetakh)&&{}-{}\sigmatr& & = 0 \\
                    &\zetakh&&{}-{}\frac{\lambda}{\sqrttt(\sigmay + \zetaih)}HP(\sigma-\zetakh)&&{}-{}\zetakhtr&& = 0
                \end{alignedat} \\
                {\psicih}'(\zetaih) - \lambda\sqrttt - \bzetaih = 0 \\
                \sqrt{ (\sigma-\zetakh)^\top P (\sigma-\zetakh) } = \sqrttt(\sigmay + \zetaih)
            \end{gathered}
            \label{eq:elastoplastic-KKT-rewritten}
        \end{equation}
        where $\bzetaih$ is the component $\bzeta$ corresponding to $\zetaih$.
        Introducing the coordinate transformation $\sigma = Q\ys$ and $\zetakh = Q\ykh$,
        subtracting \subeq{eq:elastoplastic-KKT-rewritten}{2} from \subeq{eq:elastoplastic-KKT-rewritten}{1},
        and multiplying through from the left by $Q^\top$, we have
        \begin{equation}
            \ys - \ykh = \del{I + \frac{\lambda}{\sqrttt(\sigmay + \zetaih)}\LambdaCHP}^{-1}Q^\top(\sigmatr - \zetakhtr)
            \label{eq:y-calc}
        \end{equation}
        where the diagonal matrix $\LambdaCHP = Q^\top (C+H)P Q$.
        Solving equation \subeq{eq:elastoplastic-KKT-rewritten}{3} for $\lambda$ and substituting, we get
        \begin{equation*}
            \ys - \ykh = \del{I + \frac{3}{2}\frac{{\psicih}'(\zetaih) - \bzetaih}{\sigmay + \zetaih}\LambdaCHP}^{-1}Q^\top(\sigmatr - \zetakhtr)
        \end{equation*}
        Lastly, substituting this in \subeq{eq:elastoplastic-KKT-rewritten}{4} and squaring, we get the following scalar equation for $\zetaih$.
        \begin{equation}
            \sum_r (\Lambda_P)_r \del{\frac{[Q^\top(\sigmatr - \zetakhtr)]_r}{\sqrttt(\sigmay + \zetaih) + \sqrt{\frac{3}{2}}({\psicih}'(\zetaih) - \bzetaih)(\LambdaCHP)_r}}^2 = 1
            \label{eq:zetaih-scalar-eqn}
        \end{equation}
        where $r$ goes from 1 to 6 in 3D and from 1 to 3 in 2D. Equation \eqref{eq:zetaih-scalar-eqn} can also be written equivalently as
        \begin{equation}
            \sum_r (\Lambda_P)_r \del{\frac{[Q^\top(\sigmatr - \zetakhtr)]_r}{\sqrttt(\sigmay + {\psiih}'(\alphih)) + \sqrt{\frac{3}{2}}(\alphih - \bzetaih)(\LambdaCHP)_r}}^2 = 1
            \label{eq:zetaih-scalar-eqn-alt}
        \end{equation}
        where $\alphih$ is the kinematic internal state for isotropic hardening, and $\psiih$ is the Legendre transform of $\psicih$.
        $\alphih$ is the equivalent plastic strain \cite{SivaselvanNLM2013}.
        Equation \eqref{eq:zetaih-scalar-eqn-alt} is a scalar non-linear equation in $\alphih$, which can be solved by Newton's method.
        The other states can then be calculated as summarized in Procedure \ref{alg:return-map}.
\end{enumerate}
The consistent tangent can be calculated using \eqref{eq:Kbar-direct}.
The above process of solving \eqref{eq:IP-primal} for the constitutive model \eqref{eq:elastoplastic} is summarized in Procedure \ref{alg:return-map}.

We recount that a simpler strategy to solve \eqref{eq:IP-primal} for the constitutive model \eqref{eq:elastoplastic} is enabled by the following features.
\begin{enumerate}
    \item The model has only one yield condition
    \item $\sigma$ and $\zeta$ are uncoupled in the complementary stored energy function $\psic$
    \item The part of $\zeta$ that appears in a non-quadratic manner in $\psic$, namely the isotropic hardening state $\zetaih$, is a scalar.
        This helps in reducing the problem to a scalar equation as in \eqref{eq:zetaih-scalar-eqn-alt}.
    \item $C$, $H$ and $P$ share the same eigenvectors, so that there is a diagonalizing coordinate transformation $Q$.
\end{enumerate}
In addition, in the 3D case, $(\Lambda_P)_1 = 0$ since the von Mises yield constraint applies to the deviatoric part of the stress,
and $(\LambdaCHP)_{2:6} = 2(G + G_\text{kh})$ where $G$ and $G_\text{kh}$ are the elastic and kinematic hardening shear moduli.
The solution of \eqref{eq:IP-primal} can also be interpreted as a radial projection of the trial state on the elastic region,
and is referred to as radial return mapping \cite{SimoTaylor1985,SimoHughesCompInelasticity}.
In the 2D plane stress case, squaring the yield condition as done in \cite{SimoTaylorPlaneStress} further simplifies the expression for the consistent tangent.
However, here we keep the yield function as shown in \subeq{eq:elastoplastic}{2},
so that it remains a convex function and so that the 3D and plane stress cases can be treated in a uniform manner.
We also note that if isotropic hardening were also linear, then the objective function in \eqref{eq:IP-primal} can be written as
\begin{equation*}
    \half (\sigma - \sigmatr)^\top C^{-1} (\sigma - \sigmatr) + \half (\zeta - \zetatr)^\top H^{-1} (\zeta - \zetatr)
\end{equation*}
Thus the solution can be interpreted as the closest point projection in the $C^{-1}, H^{-1}$ norm of the trial state on the elastic region.
In the absence of such simplifying features, \eqref{eq:IP-primal} must be solved by a general constrained convex optimization strategy.

\nomenclature[$bzetakh$]{$\bzetakh$, $\bzetaih$}{Components of $\bzeta$ corresponding to $\zetakh$ and $\zetaih$,
                                                 used following equation \eqref{eq:elastoplastic-KKT} in section \ref{sec:return-map}}
\nomenclature{$\sigmatr$}{Elastic trial stress (equation \eqref{eq:trial-state})}
\nomenclature{$\zetatr$, $\zetakhtr$, $\zetaihtr$}{Elastic trial state (equation \eqref{eq:trial-state})}
\nomenclature{$Q$}{Diagonalizing coordinate transformation matrix (Table \ref{tab:PQ})}
\nomenclature[$y$]{$\ys$, $\ykh$}{Components of stress and kinematic hardening internal state in volumetric-deviatoric coordinates}
\nomenclature{$I$}{Identity matrix of appropriate size}
\nomenclature{$\Lambda_\square$}{Diagonal matrix containing the eigenvalues of the symmetric matrix $\square$, used in equation \eqref{eq:zetaih-scalar-eqn}}
\nomenclature{$r$}{Index used in equation \eqref{eq:zetaih-scalar-eqn}}
\nomenclature{$\alphih$}{Isotropic hardeing internal state conjugate to $\zetakh$, used in equation \eqref{eq:zetaih-scalar-eqn-alt}}
\nomenclature{$\psiih$}{Legendre transform of $\psicih$, used in equation \eqref{eq:zetaih-scalar-eqn-alt}}
\nomenclature{$G$, $G_\text{kh}$}{Elastic and kinematic hardening shear moduli used only at the end of section \ref{sec:return-map}}

\begin{algorithm}
    \begin{algorithmic}[1]
        \State Given $\sigman$, $\zetakh^n$, $\zetaih^n$, $\mu$
        \State Compute $\sigmatr = \sigman + CB\mu$, $\zetakhtr = \zetakh^n$ and $\zetaihtr = \zetaih^n$
        \If{$\varphi(\sigmatr,\zetatr) \le 0$}
            \State $\sigma = \sigmatr$, $\zeta = \zetatr$, $\lambda = 0$
            \State $\Kbar = C$
        \Else
            \State Compute $Q^\top(\sigmatr - \zetatr)$
            \State Solve \eqref{eq:zetaih-scalar-eqn-alt} for $\alphih$
            \State $\zetaih = {\psiih}'(\alphih)$
            \State Calculate $\lambda$ using \subeq{eq:elastoplastic-KKT-rewritten}{3} 
            \State Compute $\sigma$ and $\zetakh$ from equation \subeq{eq:elastoplastic-KKT-rewritten}{1,2}
            \State Calculate $\Kbar$ using equation \eqref{eq:Kbar-direct}
        \EndIf
    \end{algorithmic}
    \caption{\textsc{Solution of \eqref{eq:IP-primal} for constitutive model \eqref{eq:elastoplastic} (return map)}}
    \label{alg:return-map}
\end{algorithm}

\section{Numerical example}\label{sec:NumExample}
\begin{figure}
    \centering
    \resizebox{0.4\textwidth}{!}{\includegraphics{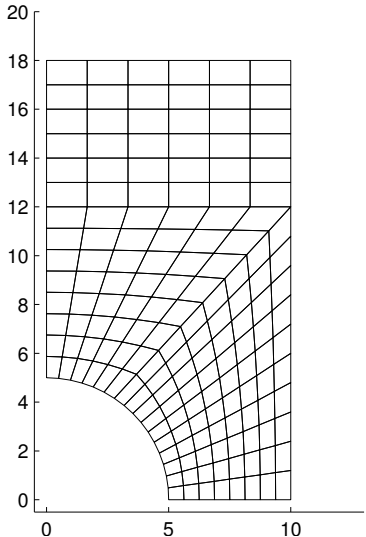}}
    \caption{Finite element model of a plate with circular hole}
    \label{fig:plate-with-hole}
\end{figure}

In this section, we present a numerical example to illustrate the application of Procedure \ref{alg:Newt}.
We use a finite element model of a plate with a circular hole subject to uniform extension at the end,
considered by Simo and Taylor \cite{SimoTaylorPlaneStress}. Due to symmetry, only a quarter of the plate needs to be represented as shown in Figure \ref{fig:plate-with-hole}.
In order to compare the numerical solution and convergence characteristics, we use the same displacement-based finite element and discretization as used in this reference.
We note however, that the formulation presented in this paper can be used with other finite element types.
A constitutive model characterized by linear elasticity, linear isotropic hardening, and a single von Mises' yield function, is used to represent material behavior.
This is a special case of the material model \eqref{eq:elastoplastic} with $\psicih(\zetaih) = \thalf h^{-1} \zetaih^2$.
Thus the Procedure \ref{alg:return-map} can be used to solve \eqref{eq:IP-primal}.
Material properties are chosen following \cite{SimoTaylorPlaneStress} as $E = 70$, $\nu=0.2$, $\sigmay=0.243$ and $h=2.24$,
where $E$ is Young's modulus and $\nu$ is Possion's ratio.
Extension is applied to the plate by imposing a uniform displacement at the top edge of the plate.
Of the three sequences of increments considered in reference \cite{SimoTaylorPlaneStress},
we only present results for the most challenging case, where two displacements increments of size 0.5 are applied in succession.
In Procedure \ref{alg:func-grad-hess}, the backtracking parameters are taken as $\beta = 10^{-4}$ and $\gamma = 0.5$.
\nomenclature{$E$}{Young's modulus (used in the numerical examples)}
\nomenclature{$\nu$}{Poisson's ratio (used in the numerical examples)}
\nomenclature{$h$}{Isotropic hardening modulus used in the numerical example}

\begin{figure}
    \centering
    \subfloat[]
    {
        \resizebox{0.4\textwidth}{!}{\includegraphics[clip]{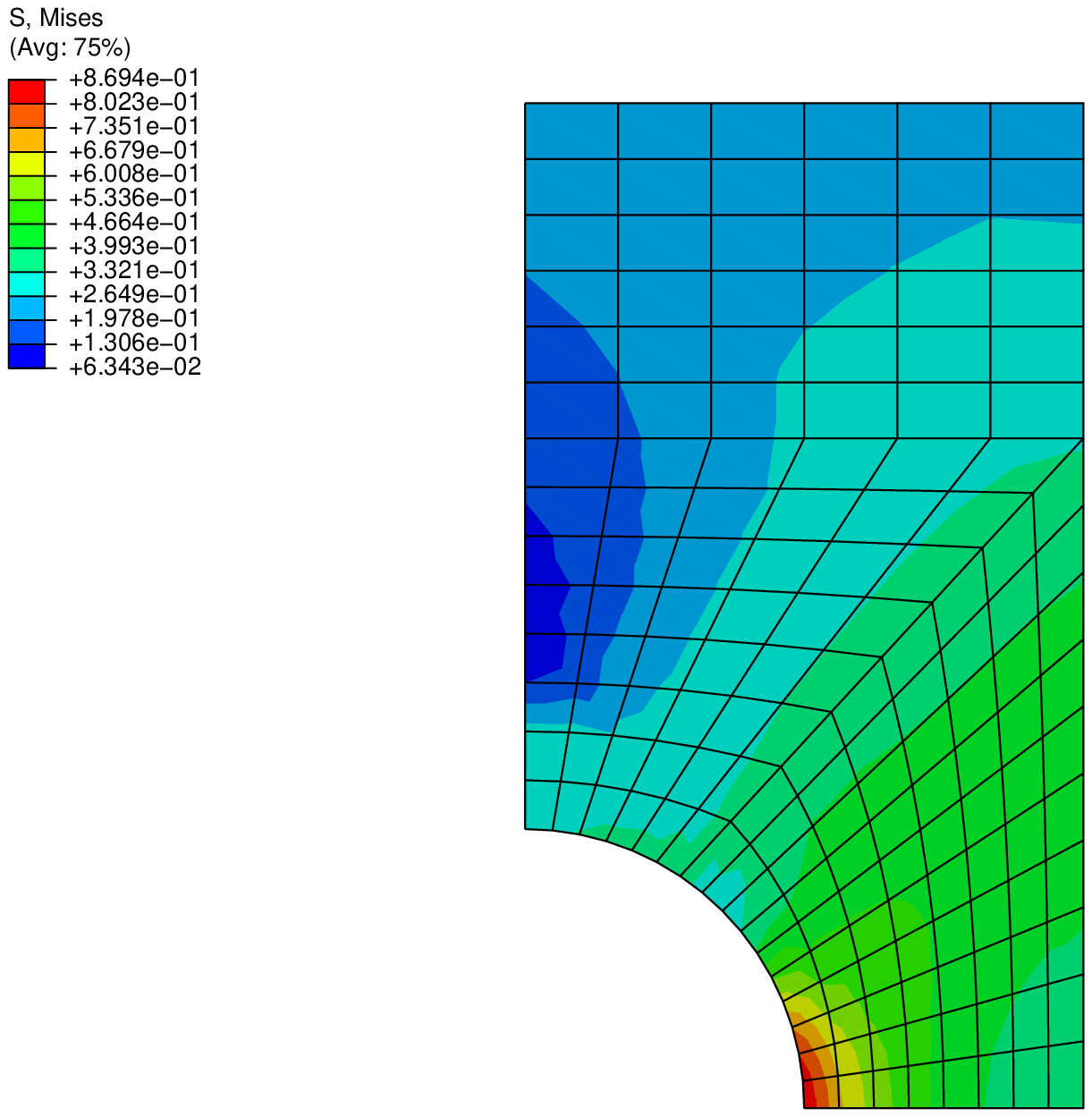}}
        \label{fig:contours:vonMises}
    }
    \hspace{0.1\textwidth}
    \subfloat[]
    {
        \resizebox{0.207\textwidth}{!}{\includegraphics[clip]{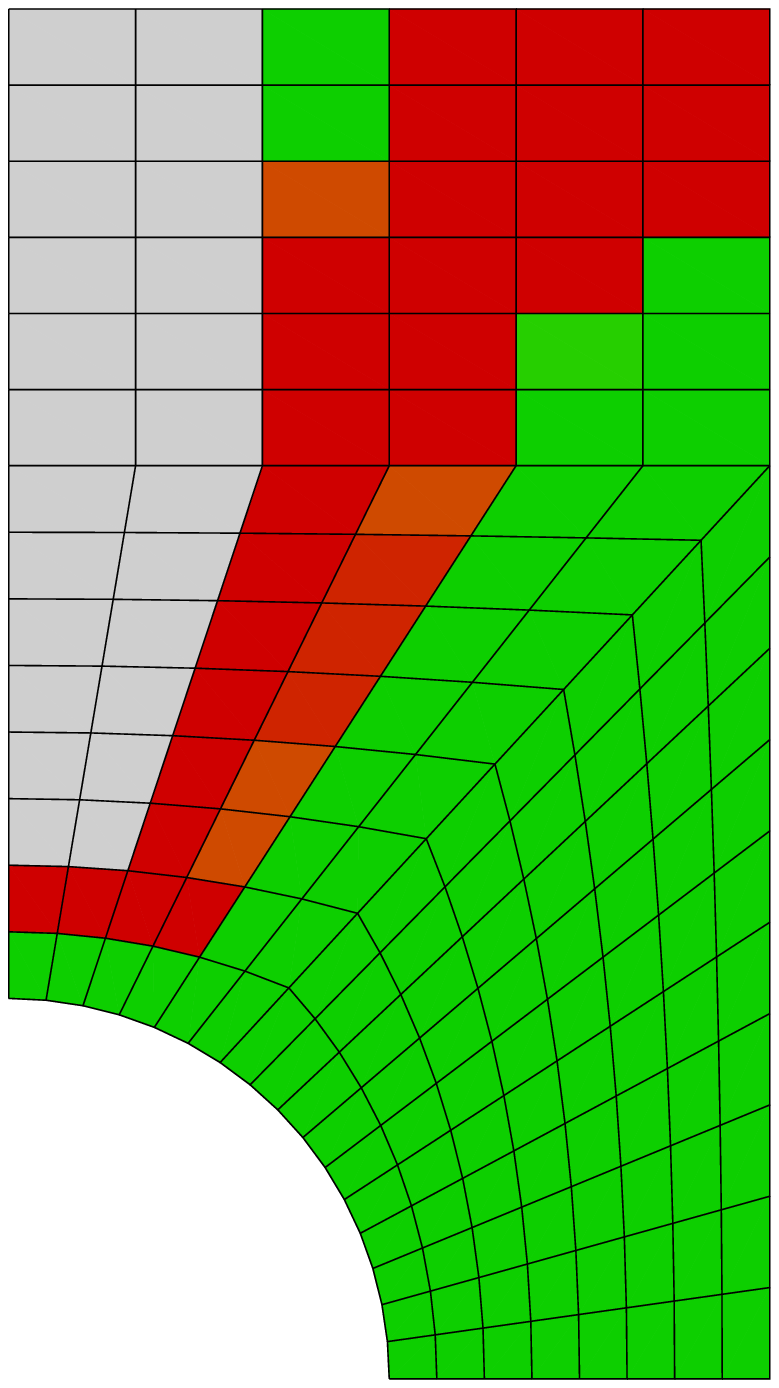}}
        \label{fig:contours:plasticity}
    }
    \caption{Numerical example (a) von Mises stress, $\sqrt{\sigma^\top P \sigma}$, after increment 2
             (b) Spread of plasticity (green --- yielded in increment 1; red --- yielded in increment 2; orange --- some integration points yielded in increment 2)}
    \label{fig:contours}
\end{figure}

Table \ref{tab:norms} shows norms of the residuals that arise in applying the optimization algorithm of Procedure \ref{alg:Newt} for each increment.
The approximate doubling of significant digits in the residuals in later iterations suggests quadratic rate of convergence.
The table also shows the step length used in each iteration, from which we see
that at most two backtracking steps used and that the full Newton step is always used in the last several iterations.
The performance is comparable to that reported in \cite{SimoTaylorPlaneStress},
with fewer backtracking steps used attributable to use of the sufficient decrease condition.
The von Mises stress computed and the spread of plasticity are shown in Figure \ref{fig:contours}.

\begin{table}
   \centering
   \begin{tabular}{rrlrl}
        \toprule
                        & \multicolumn{4}{c}{Increment (Step Length)} \\
            Iteration   & \multicolumn{2}{c}{1}	        & \multicolumn{2}{c}{2}	       	\\
        \midrule
            1	&	6.66e+00	&		    &	5.93e+00	&		    \\
            2	&	2.18e+00	&		    &	1.44e+00	&	(0.5)	\\
            3	&	1.38e+00	&	(0.25)	&	6.80e-01	&	(0.5)	\\
            4	&	1.57e+00	&	(0.5)  	&	7.00e-01	&	(0.5)	\\
            5	&	1.33e+00	&	(0.25)	&	5.33e-01	&	(0.25)	\\
            6	&	1.12e+00	&	(0.25)	&	3.43e-01	&	(0.5)	\\
            7	&	8.59e-01	&		    &	1.34e-01	&		\\
            8	&	3.96e-01	&	(0.5)   &	3.12e-02	&		\\
            9	&	1.75e-01	&		    &	9.63e-03	&		\\
            10	&	2.06e-02	&		    &	8.87e-04	&		\\
            11	&	1.53e-03	&		    &	2.45e-06	&		\\
            12	&	6.15e-06	&		    &	1.89e-11	&		\\
            13	&	1.42e-10	&		    &		        &		\\
        \bottomrule
    \end{tabular}
    \caption{Norm of of unbalanced force, $||\nabla\bar{\Pi}(\mu)||$}
    \label{tab:norms}
\end{table}

\begin{figure}
    \centering
     \subfloat[Increment 2, iteration 2]
    {
        \resizebox{0.6\textwidth}{!}{\includegraphics{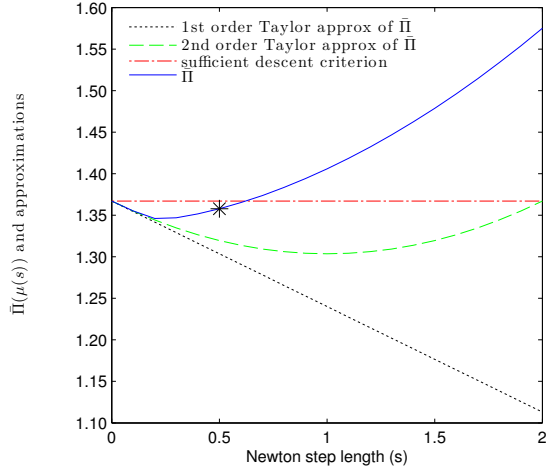}}
        \label{fig:linesearch:inc2iter2}
    } \\
    \subfloat[Increment 2, iteration 10]
    {
        \resizebox{0.6\textwidth}{!}{\includegraphics{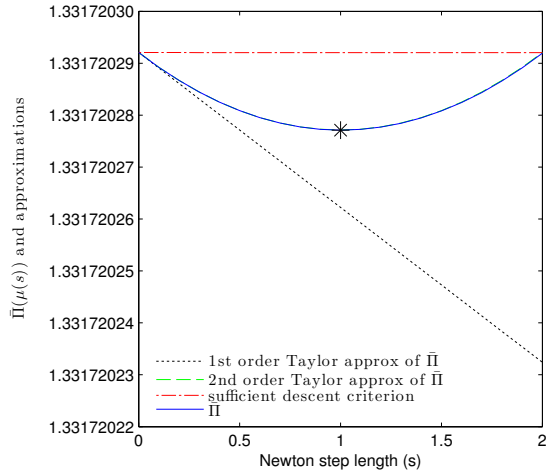}}
        \label{fig:linesearch:inc2iter10}
    }
    \caption{Visualization of backtracking for step length determination}
    \label{fig:linesearch}
\end{figure}

A typical backtracking process used to obtain the step length is depicted in Figure \ref{fig:linesearch}.
Such plots represent cross-sections along the Newton search path, and show the following as functions of the step length $s$.
\begin{enumerate}
    \item Objective function, $\bar{\Pi}(\mu + s\dmu)$
    \item First order approximation, $\bar{\Pi}(\mu) + s\nabla\bar{\Pi}(\mu)^\top\dmu$
    \item Second order approximation, $\bar{\Pi}(\mu) + s\nabla\bar{\Pi}(\mu)^\top\dmu + \thalf s^2 \dmu^\top \nabla^2\bar{\Pi}(\mu) \dmu$
    \item Sufficient descent criterion, $\bar{\Pi}(\mu) + \beta s\nabla\bar{\Pi}(\mu)^\top\dmu$
\end{enumerate}
These plots illustrate the working of the backtracking search, as well as serve to check derivative computations during implementation.
In Figure \ref{fig:linesearch:inc2iter2}, the iterate is far away from the solution.
Thus the objective function deviates considerably from its second order approximation,
and a step length of 0.5 is required to satisfy the sufficient decrease condition \eqref{eq:suff-decr}.
On the other hand in Figure \ref{fig:linesearch:inc2iter10}, the iterate is close the solution,
the second order approximation closely follows the function, and the full Newton step is taken.
We note that because of the optimization reformulation, we are able to select the step length based on
sufficient decrease of the objective function, rather than on heuristics such as number of iterations.
In this example, despite the rather large displacement increments, at most two backtracking steps are needed in any iteration
as seen from Table \ref{tab:norms}.

\section{Relationships with other approaches}\label{sec:relationships}
In section \ref{sec:Opt}, we derived three equivalent optimization problems for incremental state update
--- a primal problem \eqref{eq:Primal}, the dual problem \eqref{eq:Dual}, and a reduced dual problem \eqref{eq:DualProblem-Nested}.
In the present section, we briefly review relationships these problems share with other approaches recently presented in the literature.

There has been recent work on solving the primal problem in the context of quasi-static elastoplastic models.
Krabbenhoft et al.\ \cite{Krabbenhoft-InteriorPoint-2007} arrive at the primal problem through a different route,
and solve it directly by an interior point method.
Bilotta et al.\ \cite{Bilotta2012} also solve the primal problem directly, but by a sequential quadratic programming (SQP) approach.
Wieners \cite{Wieners2007} solves a KKT system similar to \eqref{eq:GoverningEq-discrete-KKT} by SQP.
When $\psic$ is quadratic (corresponding to linear elastic and linear hardening models),
and when the yield function is of the form of a 2-norm (as is the case with the von Mises function in section \ref{sec:return-map}),
then \eqref{eq:GoverningEq-discrete-opt} can be recast as a second order cone program (SOCP) \cite{Yonekura2012}.
This can in turn be reformulated as a semi-definite program (SDP) \cite{PracticalOptimizationAntoniouLu},
which can be solved by efficient interior point methods that have been recently developed.
Krabbenhoft et al.\ \cite{KrabbenhoftNonassociated,KrabbenhoftCamClay} also use this approach to approximate a non-associated flow model by a convex one.
Sivaselvan et al.\ \cite{Sivaselvan2009} solve a primal-type problem that arises in dynamics problems.
Further relationships to dynamic problems are discussed in \cite{SivaselvanNLM2013}.
In all of these primal approaches, a concept similar to the Sherman-Morrison-Woodbury formula is used to reduce the linear system to be solved \cite{Sivaselvan2009}.
Bilotta et al.\ \cite{Bilotta2012} discuss the relationship between the primal optimization problem
and the displacement-based nested approach through the dual problem, albeit lightly differently than done in this paper.
In \cite{Tangaramvong2012}, the KKT conditions are solved directly as a mixed complementarity problem (MCP).

As noted in section \ref{sec:reformulate:derivatives}, equation \eqref{eq:equilibrium} can also be obtained directly from \eqref{eq:GoverningEq-discrete-KKT},
since \subeq{eq:GoverningEq-discrete-KKT}{1,2,4} are the KKT conditions of \eqref{eq:IP-primal}.
However, when viewed in this manner, the minimization structure at the global level is not apparent.
In the context of plasticity models, $(\sigtildem(\mu),\zetatildem(\mu))$ is referred to as a \emph{return map}.
As discussed at the end of section \ref{sec:return-map}, for certain classes of problems, the return map defined by $\eqref{eq:IP-primal}$
can be viewed as a \emph{projection} \cite{SimoHughesCompInelasticity,Wieners2007}.
Other approaches are also possible, for example local problem is considered as a linear complementarity problem (LCP) in \cite{Bassi2012a}.
Following the seminal paper \cite{SimoTaylor1985}, Newton's method has been successfully applied to solve $\eqref{eq:equilibrium}$ in many cases.
However, since the return map is not differentiable, convergence properties of Newton's method cannot be rigorously established.
Consequently, \emph{semismooth} methods have been presented together with convergence analyses \cite{Sauter2011}.
\begin{figure}
    \centering
    \resizebox{0.95\textwidth}{!}{\includegraphics{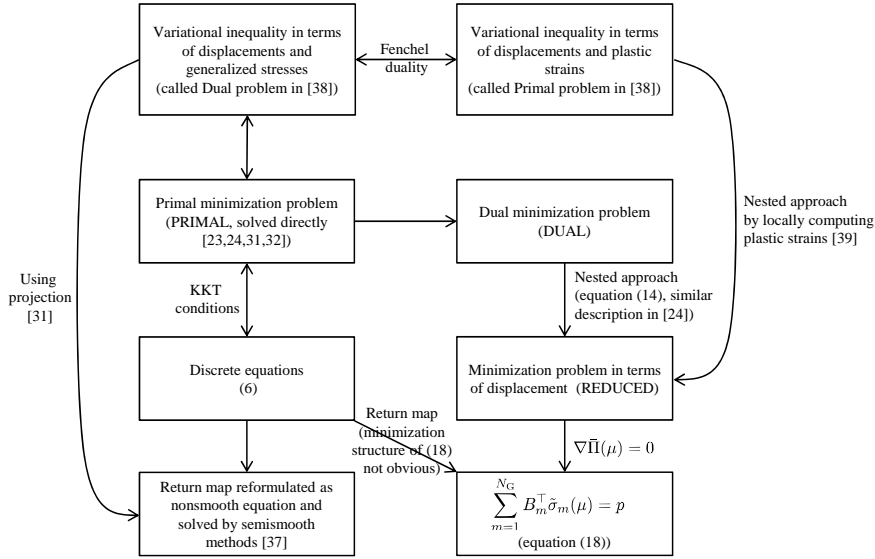}}
    \caption{Interrelation of various approaches to elastoplastic problems using optimization}
    \label{fig:interrelations}
\end{figure}

Optimization problems for quasi-static elastoplasticity can also be derived starting
from variational inequality formulations presented for example by Han and Reddy \cite{HanReddyBook}.
An optimization problem developed in this manner in terms of displacements $\mu$ and kinematic internal states or plastic strains $\alpha$
is considered by Alberty et al. \cite{Alberty1999}.
They use the optimization structure to develop error measures and an adaptive mesh-refining strategy.
Gruber and Valdman \cite{GruberSlantNewton2009} derive a minimization problem of the form
\eqref{eq:DualProblem-Nested} for a class of elastoplastic models again starting from a variational inequality.
They prove the differentiability of a function similar to $\bar{\Pi}$ using the Moreau-Yosida theorem.
Recognizing that the second derivative is not continuous, as discussed in section \ref{sec:reformulate:derivatives},
they use a slant Newton method.
Some of the interrelations described above are shown schematically in Figure \ref{fig:interrelations}.

\section{Concluding remarks}\label{sec:Conc}
We take a mathematical programming approach to the incremental state update of non-linear mechanics models
with material behavior described by stored energy and dissipation function.
In the case where these energy functions are convex and of a specific form (equation \eqref{eq:YieldFunc}),
the classical displacement-based nested approach can be reformulated as a reduced dual optimization problem \eqref{eq:DualProblem-Nested}.
We have presented an optimization reformulation of the classical displacement-based nested approach to incremental state update.
This reformulation allows visualizing the working of the algorithm by means of geometrical constructions such as Figure \ref{fig:linesearch},
and further illustrates the unifying nature of the mathematical programming approach to state update.
Connections with related algorithms recently presented in the literature have also been discussed.
We envision that such understanding gained will aid in developing algorithms for more complex non-linear material models,
such as those involving non-associated flow, softening etc., for example by means of successive optimization strategies.

\begin{thenomenclature}

 \nomgroup{Symbols}

  \item [{$B$}]\begingroup Linearized strain-displacement matrix $(\in \mathbb{R}^{\Nsig\NG\times\NDOF})$\nomeqref {1}
		\nompageref{5}
  \item [{$B_m$}]\begingroup Rows of the matrix $B$ corresponding to integration point $m$ ($\in \mathbb{R}^{\Nsig\times\NDOF}$)\nomeqref {8}
		\nompageref{9}
  \item [{$\Bprsc$}]\begingroup Deformation-displacement matrix associated with DOF with prescribed displacement\nomeqref {6}
		\nompageref{8}
  \item [{$C_d$}]\begingroup Damping matrix $(\in \mathbb{R}^{\NDOF \times \NDOF})$\nomeqref {1}
		\nompageref{5}
  \item [{$C$}]\begingroup Matrix of elastic moduli in the constitutive model equation \eqref{eq:elastoplastic}\nomeqref {25}
		\nompageref{19}
  \item [{$\mathcal{C}$}]\begingroup Convex elastic region\nomeqref {3}
		\nompageref{6}
  \item [{$\Dif$, $\Dif_q$}]\begingroup Derivative of a differentiable function; subscript $q$ denotes with respect to argument $q$\nomeqref {1}
		\nompageref{5}
  \item [{$E$}]\begingroup Young's modulus (used in the numerical examples)\nomeqref {31}
		\nompageref{23}
  \item [{$\mathcal{E}$}]\begingroup Nonlinear deformation-displacement map used only in the remarks following equation \eqref{eq:GoverningEqGeneral}\nomeqref {1}
		\nompageref{5}
  \item [{$G$, $G_\text{kh}$}]\begingroup Elastic and kinematic hardening shear moduli used only at the end of section \ref{sec:return-map}\nomeqref {31}
		\nompageref{22}
  \item [{$H$}]\begingroup Matrix of kinematic hardening moduli in the constitutive model equation \eqref{eq:elastoplastic}\nomeqref {25}
		\nompageref{19}
  \item [{$I$}]\begingroup Identity matrix of appropriate size\nomeqref {31}
		\nompageref{22}
  \item [{$\calI_m$}]\begingroup Index set of activated yield conditions defined in equation \eqref{eq:index-set}\nomeqref {20}
		\nompageref{16}
  \item [{$\Kss$, $\Ksx$, $\Kxx$}]\begingroup Derivative matrices defined in equation \eqref{eq:IP-derivative:defs}\nomeqref {20}
		\nompageref{16}
  \item [{$\Kbar_m$}]\begingroup Consistent tangent stiffness matrix for integration point $m$\nomeqref {21}
		\nompageref{16}
  \item [{$\Kbarxx$}]\begingroup Derivative matrix defined in equation \eqref{eq:Kbar-direct:defs}\nomeqref {21}
		\nompageref{16}
  \item [{$\calL$}]\begingroup Lagrangian of the primal optimization problem (equation \eqref{eq:Lagrangian})\nomeqref {8}
		\nompageref{9}
  \item [{$\calL_m$}]\begingroup Contribution to primal Lagrangian of integration point $m$ (equation \eqref{eq:LagrangianIPm})\nomeqref {10}
		\nompageref{10}
  \item [{$M$}]\begingroup Mass matrix $(\in \mathbb{R}^{\NDOF \times \NDOF})$\nomeqref {1}
		\nompageref{5}
  \item [{$\maxbacktrack$}]\begingroup Maximum number of backtracking steps in Procedure \ref{alg:Newt}\nomeqref {24}
		\nompageref{17}
  \item [{$\maxiter$}]\begingroup Maximum number of Newton iterations in Procedure \ref{alg:Newt}\nomeqref {24}
		\nompageref{17}
  \item [{$\NDOF$}]\begingroup Number of free degrees of freedom\nomeqref {1}
		\nompageref{5}
  \item [{$\NG$}]\begingroup Number of integration points\nomeqref {1}
		\nompageref{5}
  \item [{$\Nsig$}]\begingroup Number of stress components per material point (for example, 3 in 2D problems and 6 in 3D problems)\nomeqref {1}
		\nompageref{5}
  \item [{$\Nzeta$}]\begingroup Number of internal variables per material point (for example, 4 for 2D elastoplasticity with one yield condition and combined isotropic and kinematic hardening)\nomeqref {1}
		\nompageref{5}
  \item [{$P$}]\begingroup Matrix in von Mises yield condition \subeq{eq:elastoplastic}{2}\nomeqref {25}
		\nompageref{19}
  \item [{$\Phibars$, $\Phibarx$}]\begingroup Derivative matrices defined in equation \eqref{eq:Kbar-direct:defs}\nomeqref {21}
		\nompageref{16}
  \item [{$Q$}]\begingroup Diagonalizing coordinate transformation matrix (Table \ref{tab:PQ})\nomeqref {31}
		\nompageref{22}
  \item [{$\Dt$}]\begingroup Time increment\nomeqref {7}\nompageref{9}
  \item [{$\Lambda_\square$}]\begingroup Diagonal matrix containing the eigenvalues of the symmetric matrix $\square$, used in equation \eqref{eq:zetaih-scalar-eqn}\nomeqref {31}
		\nompageref{22}
  \item [{$\Phis$, $\Phix$}]\begingroup Derivative matrices defined in equation \eqref{eq:IP-derivative:defs}\nomeqref {20}
		\nompageref{16}
  \item [{$\Pi$}]\begingroup Objective function of dual problem, defined in equation \eqref{eq:DualFn}, explicit formula in equation \eqref{eq:DualFn-explicit}\nomeqref {12}
		\nompageref{11}
  \item [{$\bar{\Pi}$}]\begingroup Objective function of the reduced dual problem defined in equation \eqref{eq:Pibar}\nomeqref {16}
		\nompageref{13}
  \item [{$\Pic$}]\begingroup Objective function of the primal optimization problem (equation \eqref{eq:Pic})\nomeqref {8}
		\nompageref{9}
  \item [{$\Xi$}]\begingroup Derivative matrix defined in equation \eqref{eq:Kbar-direct:defs}\nomeqref {21}
		\nompageref{16}
  \item [{$\alpha$}]\begingroup Kinematic internal variable at a material point\nomeqref {1}
		\nompageref{5}
  \item [{$\alphih$}]\begingroup Isotropic hardeing internal state conjugate to $\zetakh$, used in equation \eqref{eq:zetaih-scalar-eqn-alt}\nomeqref {31}
		\nompageref{22}
  \item [{$\beta$}]\begingroup Sufficient decrease parameter used in choice of step length (equation \eqref{eq:suff-decr})\nomeqref {24}
		\nompageref{17}
  \item [{$\dmu$}]\begingroup Newton search direction (equation \eqref{eq:search-direction})\nomeqref {24}
		\nompageref{17}
  \item [{$\epsilon$}]\begingroup Strain at a material point\nomeqref {1}
		\nompageref{5}
  \item [{$\gamma$}]\begingroup Backtracking parameter used in Procedure \ref{alg:Newt}\nomeqref {24}
		\nompageref{17}
  \item [{$\lambda$, $\lambda_m$}]\begingroup Vector of equivalent plastic strain increments for the entire model $(\in \mathbb{R}^{\Ny\NG})$,                                       or for a single material point $(\in \mathbb{R}^{\Ny})$\nomeqref {7}
		\nompageref{9}
  \item [{$\lambdabar$, $\lambdabarnp$}]\begingroup Vector in equations \eqref{eq:phic-subgradient} and \eqref{eq:GoverningEqLagMult} $(\in \mathbb{R}^{\Ny})$;                                             corresponds to equivalent plastic strain rate; superscript denotes time increment index\nomeqref {3}
		\nompageref{6}
  \item [{$\lamtilm$}]\begingroup Minimizer of \eqref{eq:IP-dual}\nomeqref {16}
		\nompageref{13}
  \item [{$\mu$}]\begingroup Incremental displacement, $\unp-\un$\nomeqref {7}
		\nompageref{9}
  \item [{$\mu^\iter$}]\begingroup $\mu$ in iteration number $\iter$ of Newton's method (Procedure \ref{alg:Newt})\nomeqref {24}
		\nompageref{17}
  \item [{$\muprsc$}]\begingroup Increment of prescribed displacement\nomeqref {6}
		\nompageref{8}
  \item [{$\nu$}]\begingroup Poisson's ratio (used in the numerical examples)\nomeqref {31}
		\nompageref{23}
  \item [{$\phi$}]\begingroup Dissipation function for a material point\nomeqref {1}
		\nompageref{5}
  \item [{$\phic$}]\begingroup Complementary dissipation function for a material point, or for entire model, depending on context\nomeqref {1}
		\nompageref{5}
  \item [{$\varphi$}]\begingroup Yield function\nomeqref {3}
		\nompageref{6}
  \item [{$\psi$}]\begingroup Stored energy function for a material point\nomeqref {1}
		\nompageref{5}
  \item [{$\psic$}]\begingroup Complementary stored energy function for a material point, or for entire model, depending on context\nomeqref {1}
		\nompageref{5}
  \item [{$\psicih$}]\begingroup Part of complementary stored energy representing kinematic hardening in the constitutive model equation \eqref{eq:elastoplastic}\nomeqref {25}
		\nompageref{19}
  \item [{$\psiih$}]\begingroup Legendre transform of $\psicih$, used in equation \eqref{eq:zetaih-scalar-eqn-alt}\nomeqref {31}
		\nompageref{22}
  \item [{$\sigma$, $\sigman$, $\sigmanp$, $\sigmam$}]\begingroup Stress at a material point $(\in \mathbb{R}^{\Nsig})$                                                           or collection over all material points $(\in \mathbb{R}^{\Nsig\NG})$ depending on context;                                                           superscripts denote time increment index, and subscripts material point index\nomeqref {1}
		\nompageref{5}
  \item [{$\sigtildem$}]\begingroup Minimizer of \eqref{eq:IP-primal}\nomeqref {16}
		\nompageref{13}
  \item [{$\sigmatr$}]\begingroup Elastic trial stress (equation \eqref{eq:trial-state})\nomeqref {31}
		\nompageref{22}
  \item [{$\sigmay$}]\begingroup Uniaxial yield stress in von Mises yield condition \subeq{eq:elastoplastic}{2}\nomeqref {25}
		\nompageref{19}
  \item [{$\sigmstar$}]\begingroup Minimizer of equation \eqref{eq:DualFn} defined implicity in equations \eqref{eq:Def-SigZetaStar}\nomeqref {12}
		\nompageref{11}
  \item [{$\zeta$, $\zetan$, $\zetanp$, $\zetam$}]\begingroup Generalized stress internal variable at a material point $(\in \mathbb{R}^{\Nzeta})$                                               or collection over all material points $(\in \mathbb{R}^{\Nzeta\NG})$ depending on context;                                               superscripts denote time increment index, and subscripts material point index\nomeqref {1}
		\nompageref{5}
  \item [{$\zetakh$, $\zetaih$}]\begingroup Components of internal variable $\zeta$ corresponding to kinematic and isotropic hardening, used in constitutive model equation \eqref{eq:elastoplastic}\nomeqref {25}
		\nompageref{19}
  \item [{$\zetamstar$}]\begingroup Minimizer of equation \eqref{eq:DualFn} defined implicity in equations \eqref{eq:Def-SigZetaStar}\nomeqref {12}
		\nompageref{11}
  \item [{$\zetatildem$}]\begingroup Minimizer of \eqref{eq:IP-primal}\nomeqref {16}
		\nompageref{13}
  \item [{$\zetatr$, $\zetakhtr$, $\zetaihtr$}]\begingroup Elastic trial state (equation \eqref{eq:trial-state})\nomeqref {31}
		\nompageref{22}
  \item [{$\bsigma$, $\bsigmam$}]\begingroup Defined in equation \eqref{eq:GoverningEq-discrete-KKT}; components of $\bsigma$ corresponding to integration point $m$\nomeqref {7}
		\nompageref{9}
  \item [{$\bzeta$, $\bzetam$}]\begingroup Defined in equation \eqref{eq:GoverningEq-discrete-KKT}; components of $\bzeta$ corresponding to integration point $m$\nomeqref {7}
		\nompageref{9}
  \item [{$\bzetakh$, $\bzetaih$}]\begingroup Components of $\bzeta$ corresponding to $\zetakh$ and $\zetaih$,                                                  used following equation \eqref{eq:elastoplastic-KKT} in section \ref{sec:return-map}\nomeqref {31}
		\nompageref{22}
  \item [{$h$}]\begingroup Isotropic hardening modulus used in the numerical example\nomeqref {31}
		\nompageref{23}
  \item [{$\iter$}]\begingroup Iteration count in Newton's method in Procedure \ref{alg:Newt}\nomeqref {24}
		\nompageref{17}
  \item [{$k$}]\begingroup Yield function component index within an integration point $\in \{1,\ldots,\Ny\}$\nomeqref {20}
		\nompageref{16}
  \item [{$m$}]\begingroup Integration point index $\in \{1,\ldots,\NG\}$\nomeqref {8}
		\nompageref{9}
  \item [{$n$}]\begingroup Time increment index\nomeqref {7}
		\nompageref{9}
  \item [{$\nback$}]\begingroup Number of backtracking steps in Procedure \ref{alg:Newt}\nomeqref {24}
		\nompageref{17}
  \item [{$p$, $\pnp$}]\begingroup External load vector, and its values at time $n+1$ $(\in \mathbb{R}^{\NDOF})$\nomeqref {1}
		\nompageref{5}
  \item [{$r$}]\begingroup Index used in equation \eqref{eq:zetaih-scalar-eqn}\nomeqref {31}
		\nompageref{22}
  \item [{$s$}]\begingroup Length of Newton step in \ref{alg:Newt}\nomeqref {24}
		\nompageref{17}
  \item [{$t$}]\begingroup Time\nomeqref {1}\nompageref{5}
  \item [{$\mathtt{tol}$}]\begingroup Convergence tolerances for Newton's method in Procedure \ref{alg:Newt}\nomeqref {24}
		\nompageref{17}
  \item [{$u$, $\un$, $\unp$}]\begingroup Displacement at free DOF, and its values at times $n$ and $n+1$ $(\in \mathbb{R}^{\NDOF})$\nomeqref {1}
		\nompageref{5}
  \item [{$v$, $\vn$, $\vnp$}]\begingroup Velocities at free DOF, and its values at times $n$ and $n+1$ $(\in \mathbb{R}^{\NDOF})$\nomeqref {1}
		\nompageref{5}
  \item [{$x$}]\begingroup Generic variable used in Figure \ref{fig:IndFunc}\nomeqref {3}
		\nompageref{6}
  \item [{$\ys$, $\ykh$}]\begingroup Components of stress and kinematic hardening internal state in volumetric-deviatoric coordinates\nomeqref {31}
		\nompageref{22}
  \item [{$\nabla$, $\nabla_q$}]\begingroup Gradient of a differentiable function; subscript $q$ denotes with respect to argument $q$\nomeqref {1}
		\nompageref{5}
  \item [{$\partial$, $\partial_q$}]\begingroup Subdifferential of a nonsmooth convex function; subscript $q$ denotes with respect to argument $q$\nomeqref {1}
		\nompageref{5}
  \item [{$\sqcup_\mathcal{C}$}]\begingroup Indicator function of the convex elastic region $\mathcal{C}$\nomeqref {3}
		\nompageref{6}

\end{thenomenclature}

\bibliography{HysModelMathProg,MLFWhitePaper,ProjNewt,PlasticityOpt}

\end{document}